\UseRawInputEncoding 
\documentclass[12pt]{amsart}
\usepackage[english]{babel}
\usepackage{amsmath}
\usepackage{amsthm}
\usepackage{amssymb}
\usepackage{mathrsfs}
\usepackage{enumerate}
\usepackage[notcite, final, notref]{showkeys}
\usepackage{dsfont}
\usepackage[dvips]{color}
\usepackage{xcolor}
\usepackage{xmpmulti}
\usepackage{fancyref}
\usepackage{tabularx}
\usepackage{multirow}
\usepackage{stfloats}
\usepackage{caption}
\usepackage{graphicx, subfig}
\usepackage{siunitx}
\usepackage{float}
\usepackage{colortbl,dcolumn}
\usepackage{diagbox}
\newtheorem{theorem}{Theorem}[section]

\newtheorem{lemma}[theorem]{Lemma}

\theoremstyle{definition}

\newtheorem{example}[theorem]{Example}
\usepackage[ruled,linesnumbered]{algorithm2e}

\def\D{\mathbb{D}}
\def\C{\mathbb{C}}
\def\H{\mathcal{H}}
\def\e{\text{e}}
\def\i{\text{i}}
\title[]{Phase Retrieval in Hardy Space}

\author[W. Qu]{Wei Qu}
\address{Wei QU, School of Mathematical Sciences\\
 Beijing Normal University \\
 China}
\email{quwei2019@bnu.edu.cn}

\author[Xiao-Yun Sun]{Xiao-Yun Sun*}
\address{Xiao-Yun Sun, Henan Agricultural University\\
China}
\email{xysun97@163.com}
\thanks{*Corresponding author.\\
}

\author[G.-T. Deng]{Guan-Tie Deng}
\address{Guan-Tie DENG, School of Mathematical Sciences\\
Beijing Normal University\\
China}
\email{96022@bnu.edu.cn}

\begin{document}
\maketitle
\begin{abstract}
This paper concerns the study of reconstructing a function $f$ in the Hardy space of the unit disc $\D$ from intensity measurements $|f(z)|,\ z\in \D.$ It's known as the problem of phase retrieval. We transform it into solving the corresponding outer and inner function through the Nevanlinna factorization Theorem. The outer function will be established based on the mechanical quadrature method, while we use two different ways to find out the zero points of Blashcke product, thereby computing the inner function under the assumption that the singular inner function part is trivial. Then the concrete algorithms and illustrative experiments follow. Finally, we give a sparse representation of $f$ by introducing the unwinding adaptive Fourier decomposition.
\end{abstract}

\bigskip


\bigskip

\date{today}
\tableofcontents

\section{Introduction}
The problem of phase retrieval is to determine a function $f$ from its magnitude measurements $|f|.$ Important applications include X-ray crystallography, transmission electron microscopy and coherent diffractive imaging, etc. For dealing with speech recognition and noise reduction, in \cite{BCE}, the authors construct real frames and complex frames for Hilbert spaces to realize signal reconstruction from the absolute value of the frame coefficients. With application in audio processing, in 2019, Alaifari et al. present a new paradigm for stable phase retrieval by reconstructing a signal in a Hilbert space up to a phase factor(\cite{ADGY}).
Meeting statistical learning theory, in \cite{BR},  Bahmani and Romberg propose a flexible convex relaxation for the phase retrieval problem that operates in the natural domain of the signal. They avoid the prohibitive computational cost and compete with recently developed non-convex techniques for phase retrieval. Inspired by the problem of atmospheric turbulence and worked for general objects, even with noisy Fourier modulus data, Fienup give a digital method for obtaining high-resolution imagery from interferometer data, which is a phase retrieval problem of optical-coherence theory (\cite{F}). Hennelly and Sheridan bring up method of phase retrieval using the fractional Fourier transform for the encryption of image information (\cite{HS}). In \cite{LZ}, the authors establish the condition on the Blaschke product such that the orthogonality of the function can be determined, up to a unimodular scalar, by the intensity measurements, which is a classical application in coherent diffraction imaging.  In \cite{WZGAC}, the authors develop a new algorithm called sparse Truncated Amplitude flow from a small number of magnitude measurements, which realise a sparse representation for signal reconstruction. The closest study to this paper is \cite{PLB}. The authors in \cite{PLB} discuss the problem of phase retrieval in the Hardy space of the unit disc. They show that if the function has no singular part then it is uniquely determined by its amplitude on the unit circle and by its amplitude on a circle inside the unit disk. Follow this idea and motivated by rich applications of phase retrieval problem in engineering field, we consider phaseless signal recovery in Hardy space of the unit disk, which play a fundamental role in the theory of analytic signals and they are closely related to causal signals and systems. We will introduce two different methods with their concrete algorithms for the problem in the present paper.\\

For self containing purpose, we first introduce the following basic knowledge (\cite{DQ,G,QI}). The Hardy space of the unit disc is defined as

\begin{eqnarray*}
H^2({\D})&\triangleq&\{f:{\D}\to {\C}\ |\ f\ {\rm is\ analytic\ and}\ f(z)=\sum_{k=0}^\infty c_kz^k\ {\rm with}\ \sum_{k=0}^\infty |c_k|^2<\infty\}\nonumber \\
&=&\{f:{\D}\to {\C}\ |\ f\ {\rm is\ analytic\ and}\ \sup_{0<r<1}\int_0^{2\pi}|f(r\e^{\i t})|^2dt<\infty\}.
\end{eqnarray*}

An analytic function $f(z)$ on $\D$ is in the \emph{Nevanlinna class}, i.e $f \in N$, if
the subharmonic function $\log^+ |f(z)|$ has a harmonic majorant. It follows that $H^2 \subset N.$

\begin{theorem}\label{Nevanlinna}{\rm(\cite{G})}
Let $f(z)\in N,\ f\not\equiv0$. Then
\begin{equation}\label{factorization}
  f(z) = CB(z)O(z)S1(z)/S2(z),\ |C| = 1,
\end{equation}
where $B(z)$ is a Blaschke product, $O(z)$ is an outer function, $S1(z)$ and
$S2(z)$ are singular functions. Except for the choice of the constant $C$, $|C| = 1,$
the factorization (\ref{factorization}) is unique. Every function of the form (\ref{factorization}) is in $N$.
\end{theorem}

This study concerns only functions $f\in H^2$ for which the singular part of its associated inner function $S1(z)/S2(z)$ is identical to 1. It is shown that if $f$ has no singular part then $f$ is uniquely determined by its amplitude on the unit circle $\partial \D$ and by its amplitude on a circle inside the unit disk. \\

Set $\partial\D_{r}=\{z \in \C:|z|=r\},\ 0<r<1.$

\begin{theorem}\label{unique}{\rm(\cite{PLB})}
Every $f\in H^2$ is uniquely determined (up to a unimodular constant) by its modulus on $\partial\D$ and $\partial\D_{r}$ for any arbitrary $0 < r < 1$, i.e. by the values
$$
 |f\left(\mathrm{e}^{\mathrm{i} t}\right)| \text { and }
\left|f\left(r \mathrm{e}^{\mathrm{i} t}\right)\right|, \quad t \in[0, 2\pi).
$$
\end{theorem}

According the above knowledge, we try to determine the outer and inner function part, respectively, by the amplitude measurements.\\

 The organization of this paper is as follows. In section 2, we will first introduce the Mechanical quadrature method of Hilbert transform for solving the outer function part. Two methods for finding zeros of Blaschke product are presented subsequently. We then give two concrete algorithms in section 3 and the illustrative experiments will be shown in the fourth section. Finally, we briefly review the theory of unwinding AFD and combine it with the problem of phase retrieval to obtain a sparse representation for the signal through its absolute value.

\section{Phase retrieval by Nevanlinna factorization}
Our phase retrieval approach is based on the Nevanlinna factorization in $H^2$. We first work on the outer function through calculating the Hilbert transform by the modulus of $f$ on $\partial \D.$ In the sequel, two practical methods for finding zeros will be introduced. \\

The function in $H^2$ has the form $f=OBS,$ where $O$ is the outer function given by

\begin{equation}\label{outer}
 O(z)=C\exp \left(\frac{1}{2 \pi}\int_{0}^{2 \pi} \frac{\e^{\i t}+z}{\e^{\i t}-z} \cdot \log \left|f\left(\e^{\i t}\right)\right| {d t}\right),\ \ |C|=1,
\end{equation}
and $B$ is a Blaschke product

\begin{equation}\label{inner}
B(z)=z^{m} \prod_{\left|\alpha_{k}\right| \neq 0} \frac{-\bar{\alpha}_{k}}{\left|\alpha_{k}\right|} \frac{z-\alpha_{k}}{1-\bar{\alpha}_{k} z}, \quad m=1,2,\cdots
\end{equation}

 $\alpha_k$ 's are the zeros of $f$ in the unit disc $\D.$ We assume that the singular inner function $S$ is trivial in this context. Clearly, the decomposition is unique except for the choice of the constant $C$ (\cite{G}).

\subsection{Mechanical quadrature method of Hilbert transform for outer function}
{}\ \\

It is sufficient to know $|f|$ on the unit circle to recover $O$ up to a unitary constant C using (\ref{outer}). For $f \in H^{2}(\mathbb{D}),$ and $z=\rho {\text e}^{{\text i}\theta}, 0<\rho<1, \theta \in[0,2 \pi),$ we consider the boundary value of the outer factor of $f(z),$

\begin{equation}\label{boundary}
\begin{aligned}
&\quad\ \lim _{\rho \rightarrow 1^{-}} O(z) \\
&=\lim _{\rho \rightarrow 1^{-}} \exp \left(\frac{1}{2 \pi}\int_{0}^{2 \pi} \frac{\e^{\i t}+z}{\e^{\i t}-z} \log \left|f\left(\e^{\i t}\right)\right| d t\right) \\
&=\lim _{\rho \rightarrow 1^{-}} \exp \left(\frac{1}{{2 \pi}}\int_{0}^{2 \pi}\left(\frac{1-\rho^{2}}{1-2 \rho \cos (\theta-t)+\rho^{2}}+\i \frac{2 \rho \sin (\theta-t)}{1-2 \rho \cos (\theta-t)+\rho^{2}}\right) \log \left|f\left(\e^{\i t}\right)\right| {d t}\right) \\
&=\exp \left(\log \left|f\left(\e^{\i \theta}\right)\right|+\i \mathcal{H} \log \left|f\left(\e^{\i \theta}\right)\right|\right),
\end{aligned}
\end{equation}

where $\H$ is the circular Hilbert transform, defined as

\begin{equation}\label{Hilbert}
 (\mathcal{H} f)(t)\triangleq\frac{1}{2 \pi} {\rm P. V.} \int_{0}^{2 \pi} f(x) \cot \frac{t-x}{2} d x=\lim _{\epsilon \rightarrow 0^{+}} \frac{1}{2 \pi} \int_{\epsilon<|t-x|<2 \pi} f(x) \cot \frac{t-x}{2} dx.
\end{equation}

We will introduce the Mechanical quadrature method (MQM) for calculating Hilbert transform, which can be found in \cite{CR,D1,D2,SD}. The main idea of this method is separation of singularity, can be briefly showed as:
\[(\mathcal{H} f)(t)=\frac{1}{2 \pi}\int_{0}^{2 \pi} (f(x)-f(t))\cot \frac{t-x}{2} d x+\frac{1}{2 \pi}\int_{0}^{2 \pi}f(t)\cot \frac{t-x}{2} d x.
\]

Combing the trapezoidal sequence of rules and trigonometric interpolation, there is
$(\mathcal{H} f)(t) \simeq(\mathbf{Q} f)(t),$ where

\begin{equation}\label{QF}
  (\mathbf{Q} f)(t)=\left\{\begin{array}{ll}
\frac{1}{n} \sum\limits_{j=0}^{n-1} \cot \frac{t-x_{j}}{2} f\left(x_{j}\right)+\cot \frac{n t}{2} f(t), \quad \text { if } t \neq x_{j}, \\

\frac{1}{n} \sum\limits_{j=0, j \neq k}^{n-1} \cot \frac{x_{j}-x_{k}}{2} f\left(x_{j}\right)+\frac{2}{n} f^{\prime}\left(x_{j}\right), \quad \text{if } t=x_{j},
\end{array}\right.
\end{equation}

the nodes $x_{j}=\frac{2 j \pi}{n},\ j=0,1, \ldots, n-1.$
In \cite{CR} and \cite{D1}, the authors gave the estimation and the convergence of the
quadrature formula ($\ref{QF}$), what is
$O(e^{-2nr})$ with $n\rightarrow \infty$.
The algorithm of the MQM for Hilbert transform is as follows.

\begin{algorithm*}[H]
		\caption*{Algorithm of MQM for Hilbert transform}
		\label{alg1}
		
		\KwIn{sampling points set $x=\{x_1,\ldots,x_n\}$;
           original signal $f(x)$}%
		\KwOut{Hilbert transform $(\mathbf{Q}f)(x)$ of $f(x)$}
Initialize $(\mathbf{Q}f)(x) = 0$; $X(1) =0$\;
 \For {$j=1:n$}
 {

 {
$X(j) = X(j)+\cot(\frac{1}{2}(t_{j}-x_{j}))$\;
 }
 }
 $(\mathbf{Q}f)(x) = \frac{1}{n}f(x)X+\cot\frac{nt}{2}f(t)$\;
\end{algorithm*}

Thus, according to (\ref{boundary}), we can obtain the outer function by calculating the Hilbert transform.

\subsection{Two methods for solving inner function}
{}\ \\

\subsubsection{The minimum value method}
We assume that $f\in H^2{(\D)}$ can be analytically extended to the unit circle, i.e. $f$ has a finite number of zeros in $\D$. Otherwise, the zeros of $f$ have an accumulation point on $\D$. It follows $f\equiv 0$ by a basic fact that zeroes of analytic functions are isolated. Then $f(a)=0$ implies
\begin{equation}\label{min}
a=\arg \min _{z \in \D}|f(z)|.
\end{equation}
We here have two ways to approximately fix the number of zeros for $f.$ On one hand, we keep $|f(z)|$ with $|f(z)|\leq \varepsilon,$ where $\varepsilon$ can be set as a small number such as $0.01,\ 0.001, \cdots,$ thus we obtain all points in $\D$ satisfied $|f(z)|\leq \varepsilon$ by applying (\ref{min}). On the other hand, set $f_1(z)=f(z),$ we obtain $\alpha_1$ satisfying $f_1(\alpha_1)=0$ through (\ref{min}). Combine the outer function $O(z)$ that has been computed in the last subsection, we have the reconstructed function, denoted by $g_1$ with $g_1(z)=O(z)\frac{z-\alpha_{1}}{1-\bar{\alpha}_{1} z}.$ Repeat the above step by setting $f_{2}=f_{1} \frac{1-\bar{\alpha}_{1} z}{z-\alpha_{1}}$ and apply (\ref{min}) to $f_2$, we further obtain $\alpha_2$ satisfying $f_2(\alpha_2)=0.$ The reconstructed function $g_2=O(z)\frac{z-\alpha_{1}}{1-\bar{\alpha}_{1} z}\frac{z-\alpha_{2}}{2-\bar{\alpha}_{2} z}.$
Whether continue repeating the above step depends on the reconstruction error defined by $$\operatorname{err}=\min \{\left\|f-g_k\right\|\},$$
where $g_k=O(z)\prod_k\frac{z-\alpha_{k}}{1-\bar{\alpha}_{k} z}.$ The concrete algorithm and illustrative experiments can be found in the next two sections.

\subsubsection{The para-conjugate method}
We will use the following lemma which can be found in \cite{PLB}.

\begin{lemma}\label{para}
Let $B$ be a Blaschke product and let $0 < r < 1$ be arbitrary. Then $B$ is uniquely determined (up to a unitary constant) by its modulus on the circle $\mathbb{\D}_{r}=\{z \in \mathbb{C}:|z|=r\},$ i.e. by the function $\left|B\left(r \mathrm{e}^{\mathrm{i} t}\right)\right|, t \in[0,2\pi).$
\end{lemma}

Let $\left\{\alpha_{n}\right\}_{n=1}^{\infty}$ be the zero set of $B$ and define

\begin{equation}\label{para}
  B_{r}(z)\triangleq B(r z) \quad \text { and } \quad P_{r}(z)\triangleq B_{r}(z) B_{r}^{*}(z),
\end{equation}
where $B_{r}^{*}(z)=\overline{B_{r}(1/\overline{z})}$ is the para-conjugate function of $B_r.$ According to \cite{PLB}, we can see that $\left\{\frac{\alpha_{n}}{r}\right\}_{n=1}^{\infty}$ is zero set of $B_r$ and $P_r$ has poles at $\mathcal{P}\left(P_{r}\right)=\mathcal{P}_{1} \cup \mathcal{P}_{2}$ with

\begin{equation}\label{poles}
  \mathcal{P}_{1}=\left\{\pi_{1, n}=\frac{1}{r \bar{\alpha}_{n}}\right\}_{n=1}^{\infty} \text { and } \mathcal{P}_{2}=\left\{\pi_{2, n}=r \alpha_{n}\right\}_{n=1}^{\infty}.
\end{equation}

By (\ref{para}), we have the following relationship for all $0<r<1,$

\begin{equation}\label{Pr}
  P_{r}\left(\mathrm{e}^{\mathrm{i} t}\right)=\left|B\left(r \mathrm{e}^{\mathrm{i} t}\right)\right|^{2}, \quad t \in[0, 2\pi).
\end{equation}

and (\ref{poles}) shows that $P_r$ is analytic in $\{z\in \C : r < |z| < 1/r\}$ and
the accumulation points of its poles lie on the circles $\partial{\D}_{1/r}$ and $\partial{\D}_r$. So Laurent theorem implies that $P_r$ is uniquely determined by the known values (\ref{Pr}) on the $\partial\D$, that is,

\begin{equation}\label{laurent}
  P_r(z)=\sum_{n=-\infty}^\infty c_nz^n,\ c_n=\frac{1}{2\pi}\int_0^{2\pi}P_r(\mathrm{e}^{\mathrm{i} t})\mathrm{e}^{-\mathrm{i} n t}dt
\end{equation}

Besides, from (\ref{para}), we have
$$
P_r(z)=\prod_n \frac{rz-\alpha_n}{1-{\overline{\alpha}}_nrz}\frac{1-{\overline{\alpha}}_nz/r}{z/r-\alpha_n}.
$$
Consequently, the zeros of $B$ can be determined by $\alpha_{n}=\pi_{2, n} / r$ from all poles of $P_r$ which satisfy $\left|\pi_{2, n}\right|<r$. Forming the corresponding Blaschke product

 $$B(z)=z^{m} \prod_{n=1}^{\infty} \frac{-\overline{\alpha}_{n}}{\left|\alpha_{n}\right|} \frac{z-\alpha_{n}}{1-{\overline{\alpha}}_{n} z}$$

 recovers $B$, wherein $m$ is the number of poles with $\left|\pi_{2, n}\right|=0$.

\section{Two algorithms}
\subsection{Algoritm 1}
{}\ \\

By using the Mechanical quadrature method for outer function, we are to find out the zero points of inner function through determining the minimum absolute value. This method is short for MQMV.\\

\begin{tabular}{l}
\hline
{\bfseries Algorithm 1:} MQMV\\
\hline
\;\textbf{Input}:\quad $\ \ |f|$\\
\;\textbf{Output}:\quad $f$\\
\;Initialize:\quad $k=0$\\
\;step 1:\quad\quad Compute the Hilbert transform $\mathcal{H}$ of $F(\e^{\i t})=\log |f(\e^{\i t})|$\\
\;step 2:\quad\quad Obtain outer function of the circle $O(\e^{\i t})=\e^{(F+\i\mathcal{H}F)(\e^{\i t})}$\\
\;step 3:\quad\quad Compute the value of inner function of the circle $|B(\e^{\i t})|=|f(\e^{\i t})|/|O(\e^{\i t})|$  \\
\;step 4:\quad\quad Obtain $|B(r\e^{\i t})|$ through Cauchy integral\\
\;step 5:\quad\quad \textbf{While} $|B(r\e^{\i t})|\leq 0.001$ \textbf{do}\\
\;\quad\quad\quad\quad\quad\quad\quad\quad\  k=k+1\\
\;\quad\quad\quad\quad\quad\quad\quad\quad\  $\alpha_k \leftarrow \arg\ \min\{|B(r\e^{\i t})|\}$\\
\;\quad\quad\quad\quad\quad{\bfseries end}\\
\;step 6:\quad\quad $f=OB$\\
\hline
\end{tabular}
\bigskip
\subsection{Algoritm 2}
{}\ \\

In the second algorithm, we still compute the outer function by the Mechanical quadrature method. Then by constructing a new function $P_r(z),$ we can compute the zero points of the corresponding Blaschke product. The para-conjugate method for computing the inner function, which is abbreviated by MQPC, can be understood in the following process. \\

\begin{tabular}{l}
\hline
{\bfseries Algorithm 2:} MQPC\\
\hline
\;\textbf{Input}:\quad $\ \ |f|$\\
\;\textbf{Output}:\quad $f$\\
\;Initialize:\quad Select $r,\ 0<r<1$\\
\;step 1:\quad\quad Compute the Hilbert transform $\mathcal{H}$ of $F(\e^{\i t})=\log |f(\e^{\i t})|$\\
\;step 2:\quad\quad Obtain outer function of the circle $O(\e^{\i t})=\e^{(F+\i\mathcal{H}F)(\e^{\i t})}$\\
\;step 3:\quad\quad Compute the value of inner function of the circle $|B(\e^{\i t})|=|f(\e^{\i t})|/|O(\e^{\i t})|$  \\
\;step 4:\quad\quad Obtain $|B(r\e^{\i t})|$ through Cauchy integral\\
\;step 5:\quad\quad Compute the Laurent series of $P_r(z)$ by $|B(r\e^{\i t})|^2$ with fixed $0<r<1$\\
\;step 6:\quad\quad Find the poles $\pi_{2,n}$ of $P_r(z)$ with $|z|<r$\\
\;step 7:\quad\quad Obtain the zeros of $B$ by $\pi_{2,n}/r$\\
\;step 8:\quad\quad $f=OB$\\
\hline
\end{tabular}

\section{Numerical experiments}
\begin{example}
$$f(z)=\frac{0.1867 z^{6}-0.00869 z^{5}}{(1-0.7842 z)(1-0.2669 z)}$$
\end{example}
\begin{figure}[H]
\begin{minipage}{0.3\linewidth}
  \includegraphics[width=2.2in,height=1.75in]{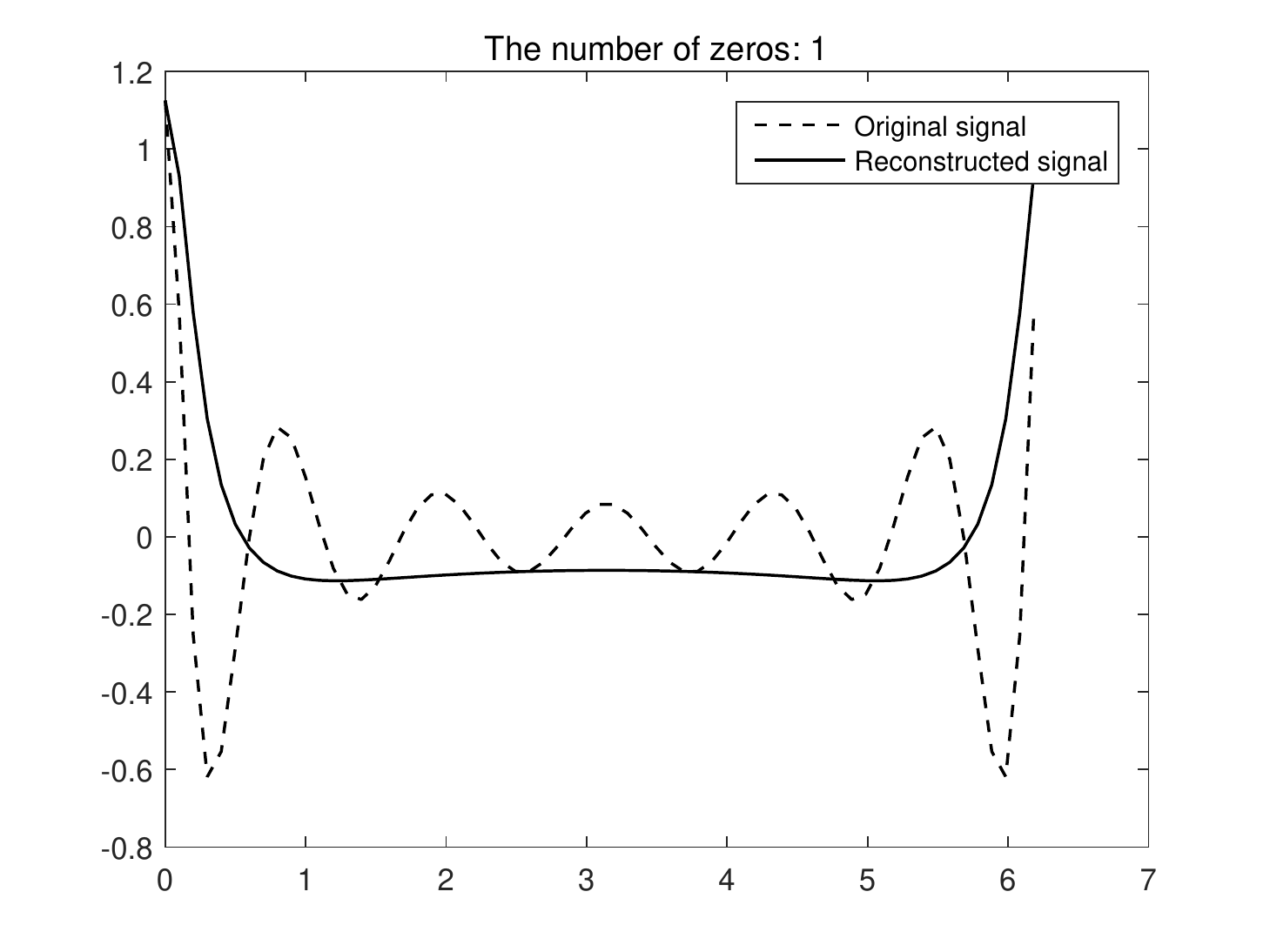}
\end{minipage}
\hfill
\begin{minipage}{0.3\linewidth}
  \includegraphics[width=2.2in,height=1.75in]{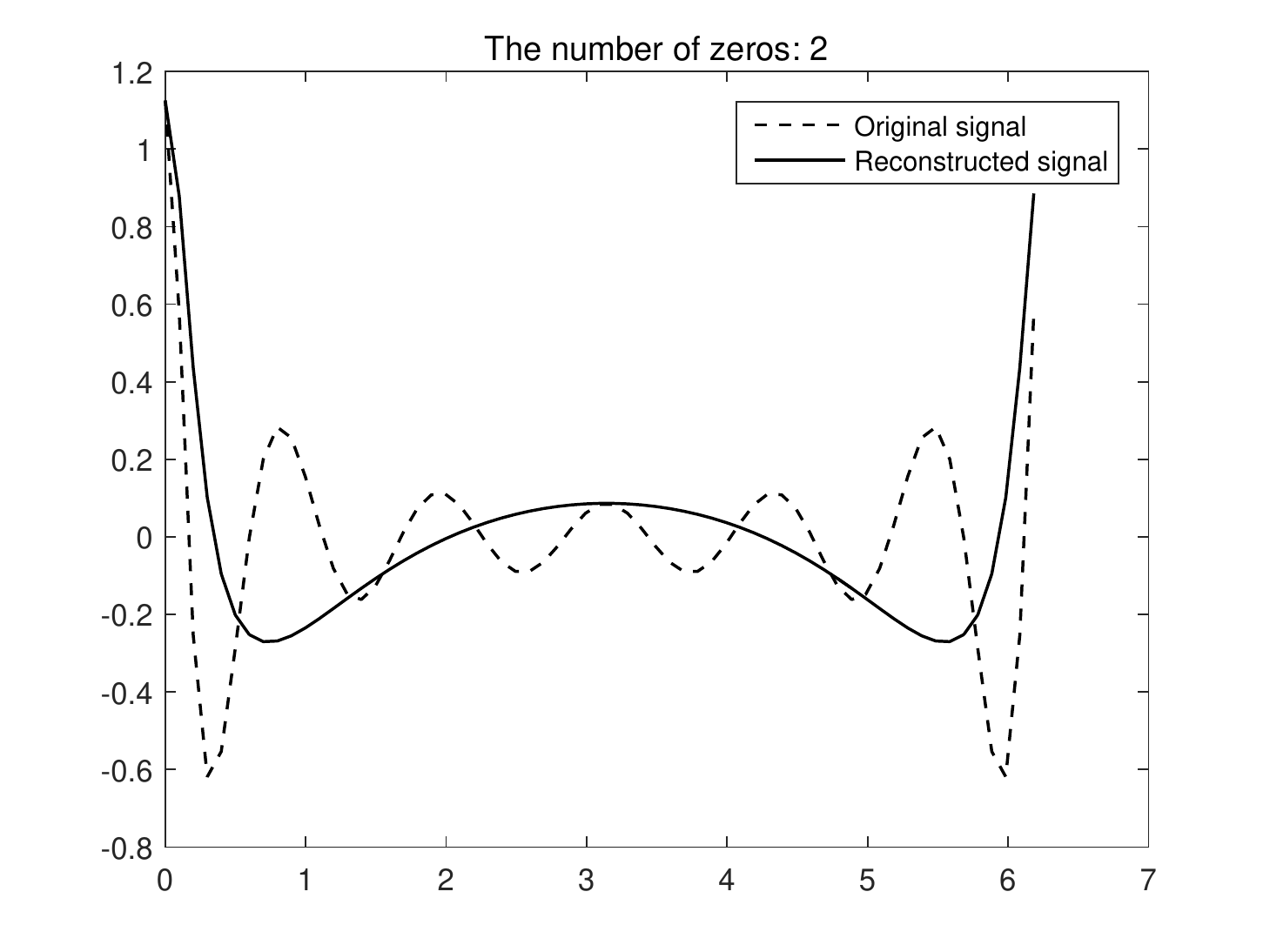}
\end{minipage}
\hfill
\begin{minipage}{0.3\linewidth}
  \includegraphics[width=2.2in,height=1.75in]{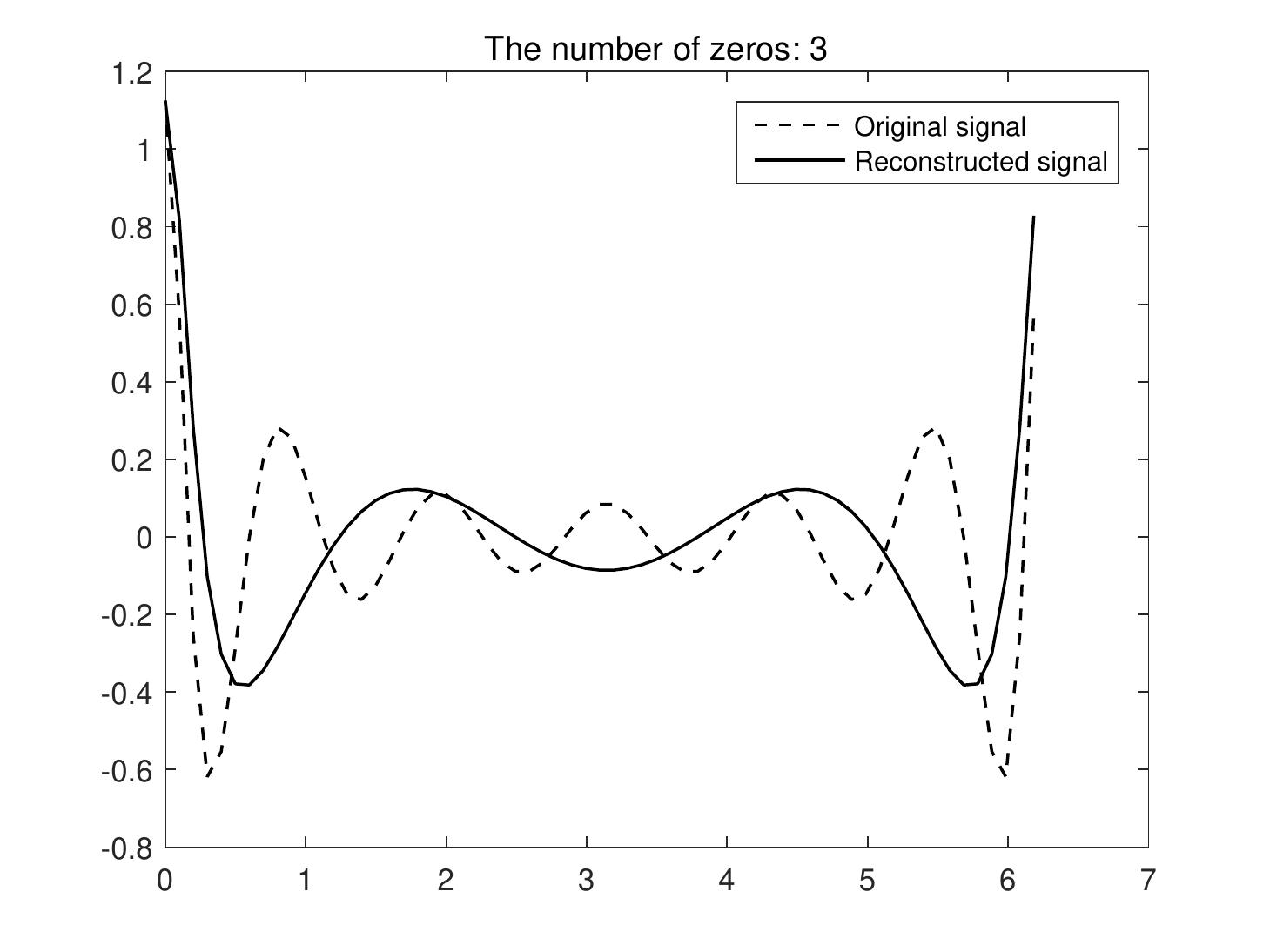}
\end{minipage}
\vfill
\begin{minipage}{0.3\linewidth}
  \includegraphics[width=2.2in,height=1.75in]{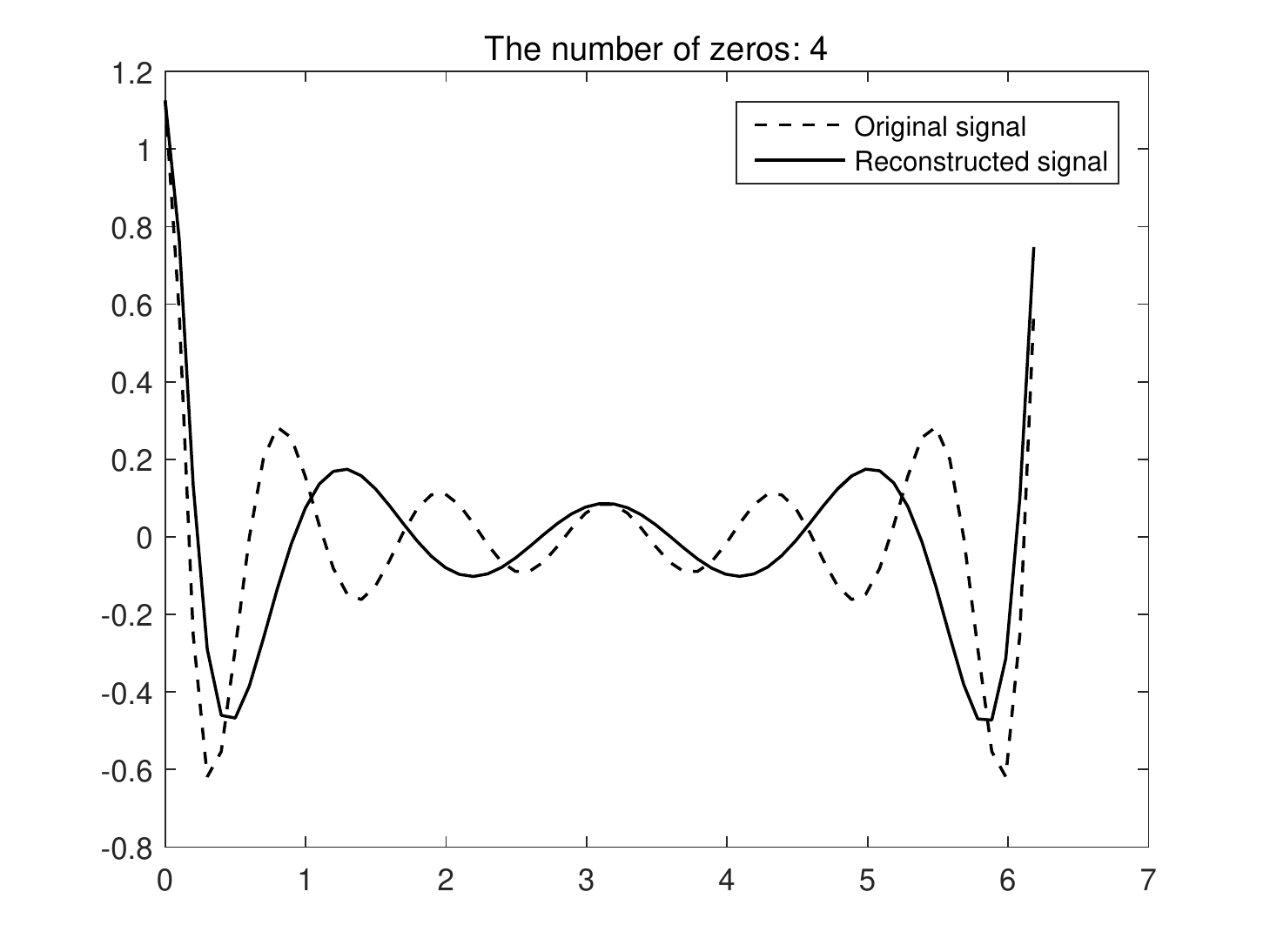}
\end{minipage}
\hfill
\begin{minipage}{0.3\linewidth}
  \includegraphics[width=2.2in,height=1.75in]{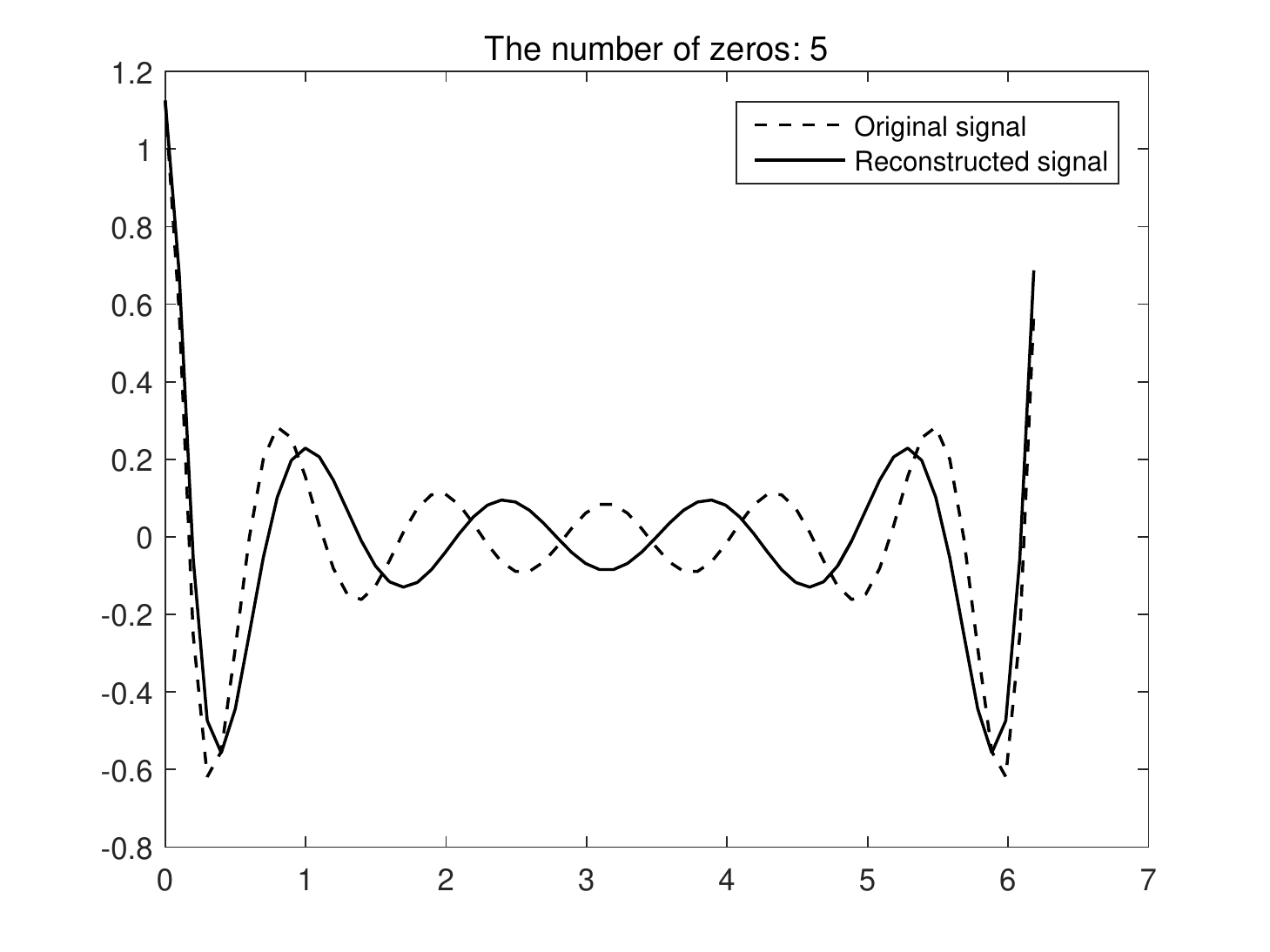}
\end{minipage}
\hfill
\begin{minipage}{0.3\linewidth}
  \includegraphics[width=2.2in,height=1.75in]{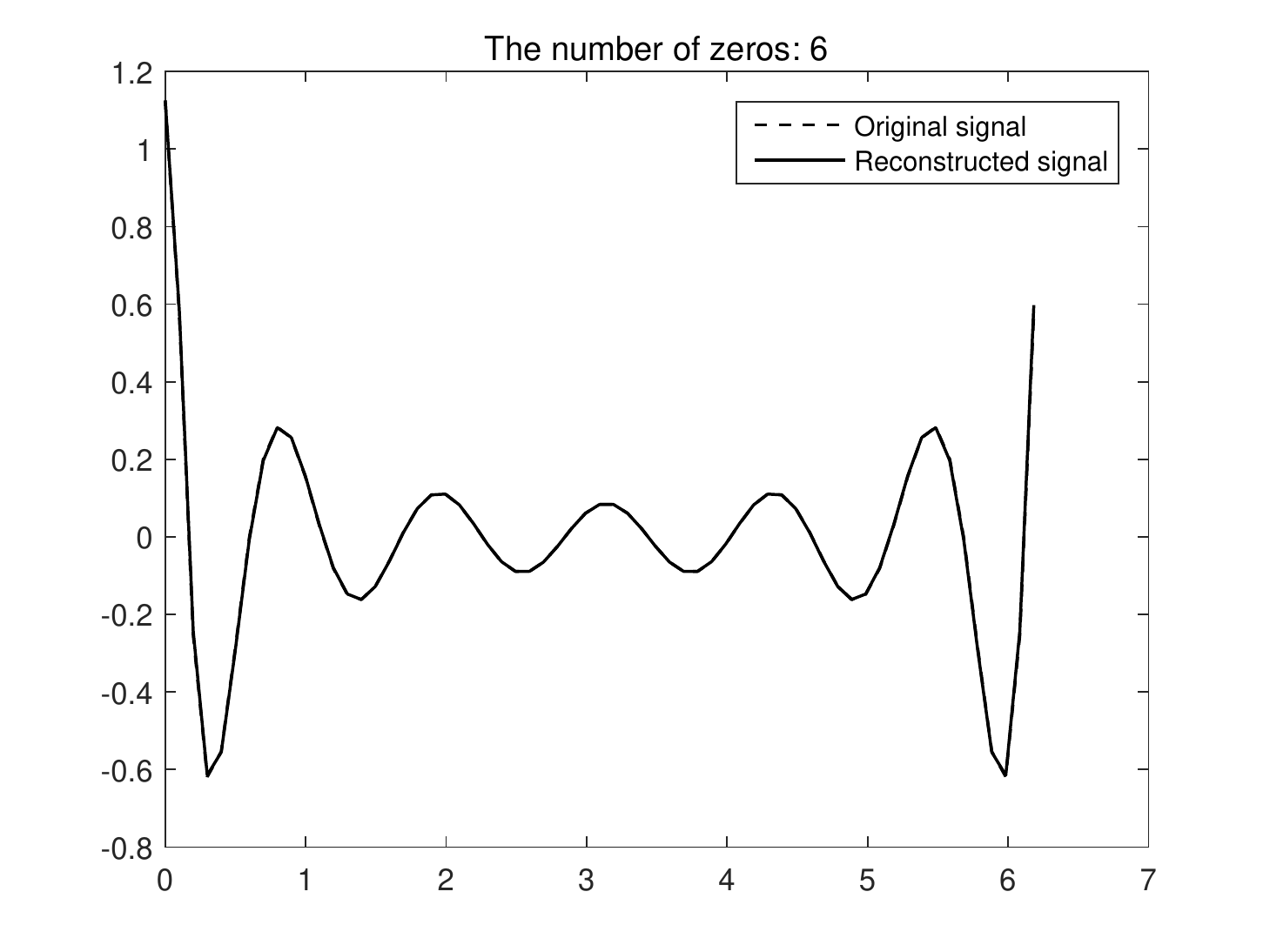}
\end{minipage}
\caption{64 sampling points}
\label{64}
\end{figure}

\begin{figure}[H]
\begin{minipage}{0.3\linewidth}
  \includegraphics[width=2.2in,height=1.75in]{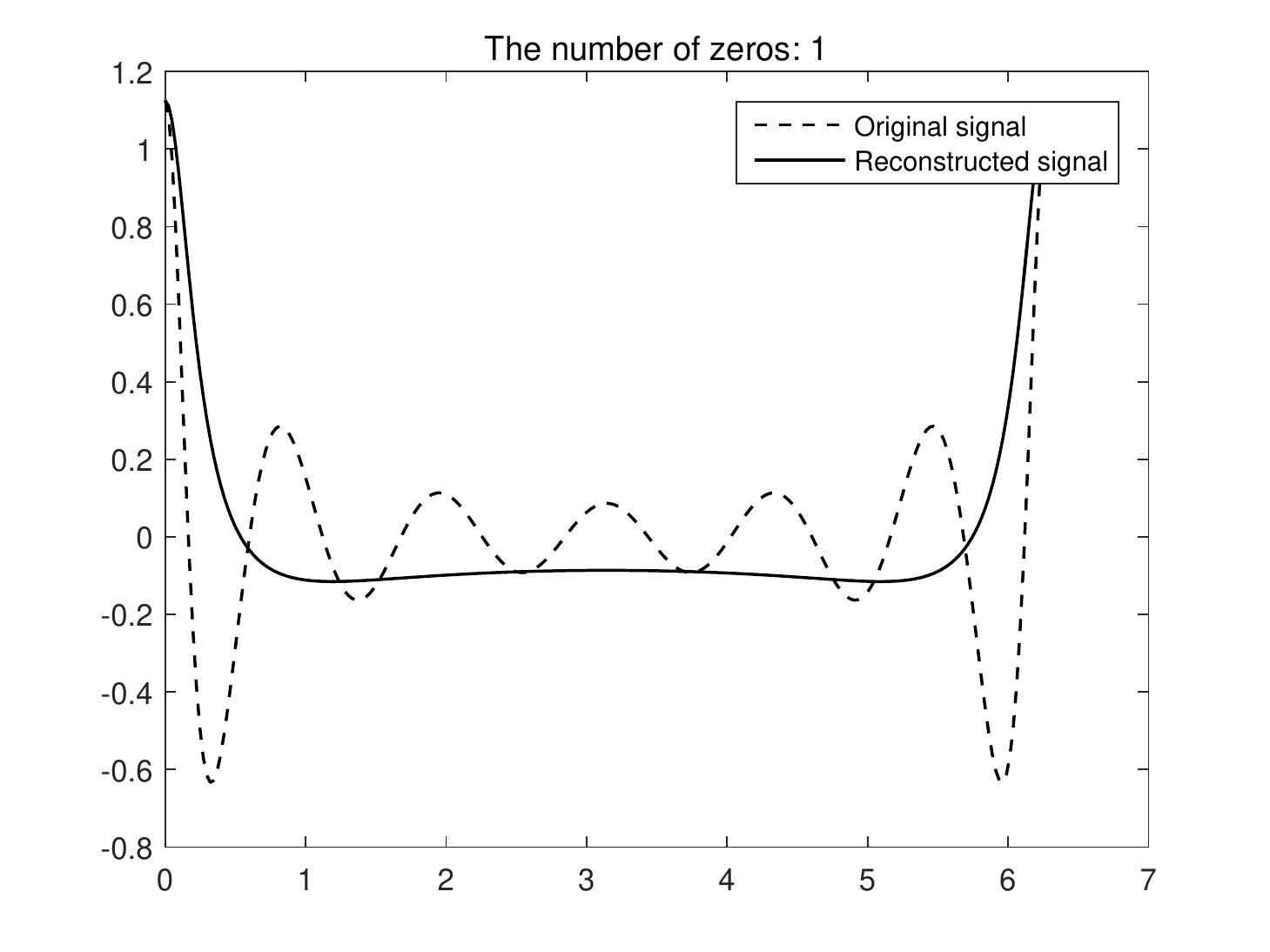}
\end{minipage}
\hfill
\begin{minipage}{0.3\linewidth}
  \includegraphics[width=2.2in,height=1.75in]{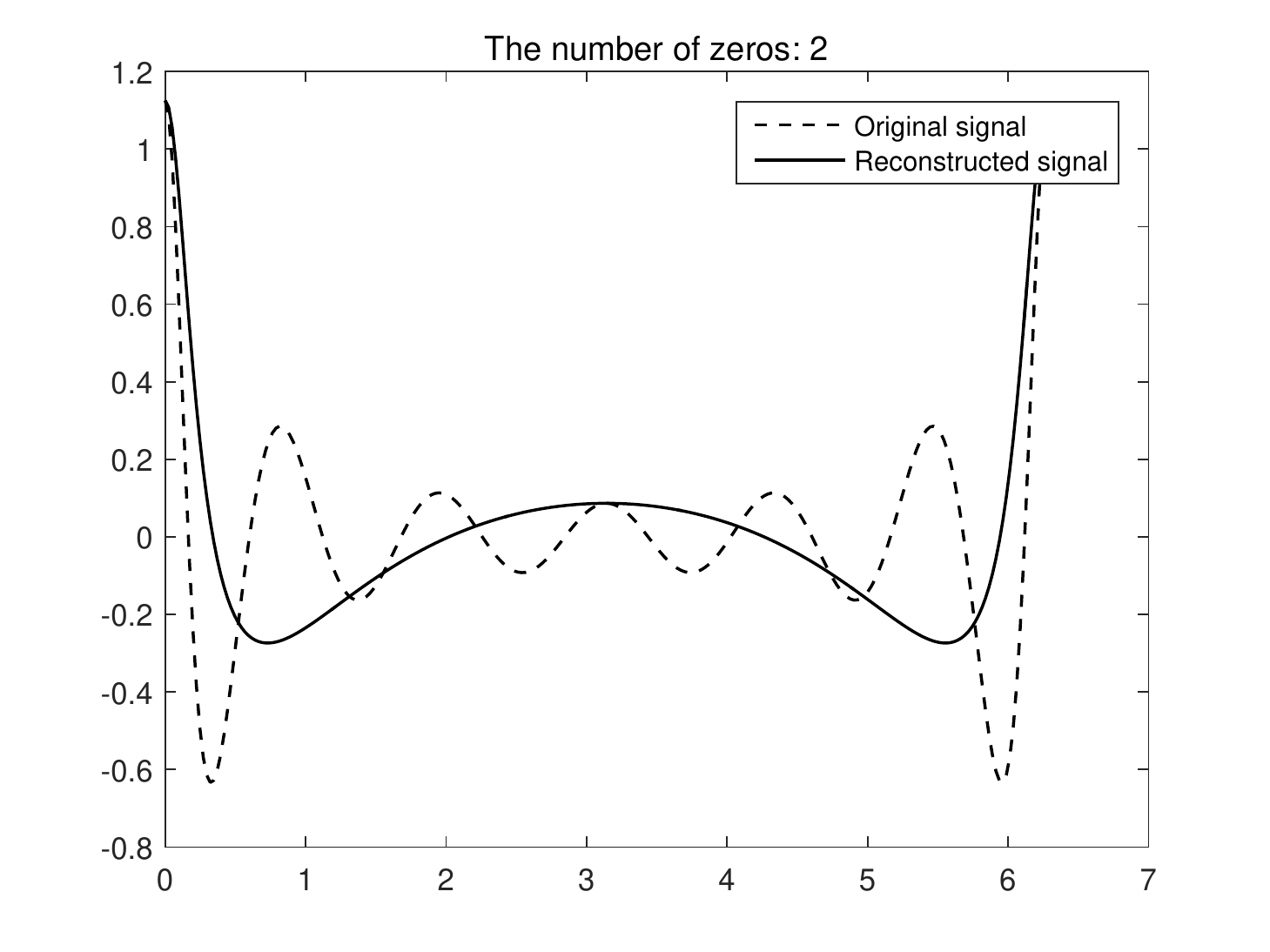}
\end{minipage}
\hfill
\begin{minipage}{0.3\linewidth}
  \includegraphics[width=2.2in,height=1.75in]{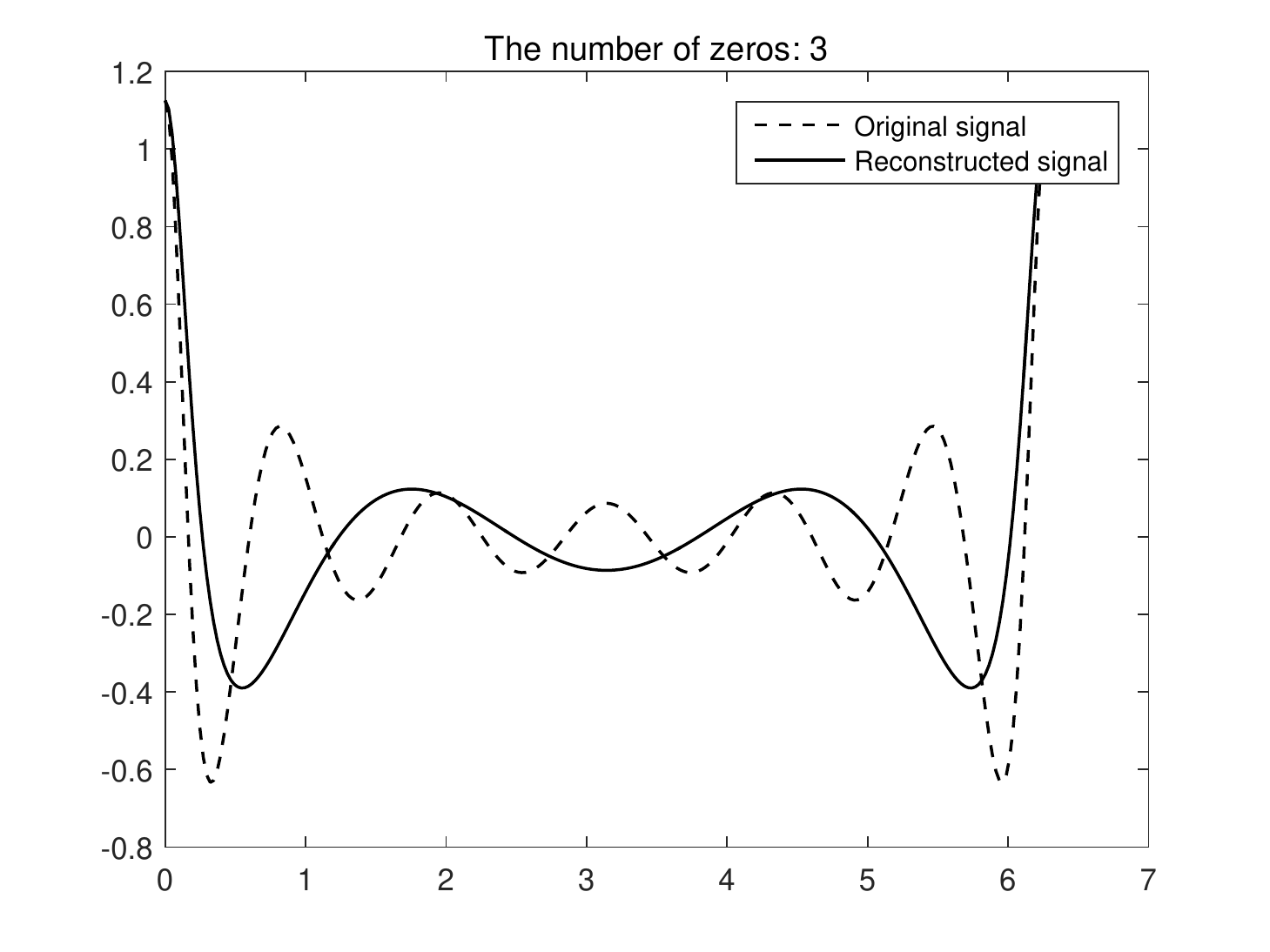}
\end{minipage}
\vfill
\begin{minipage}{0.3\linewidth}
  \includegraphics[width=2.2in,height=1.75in]{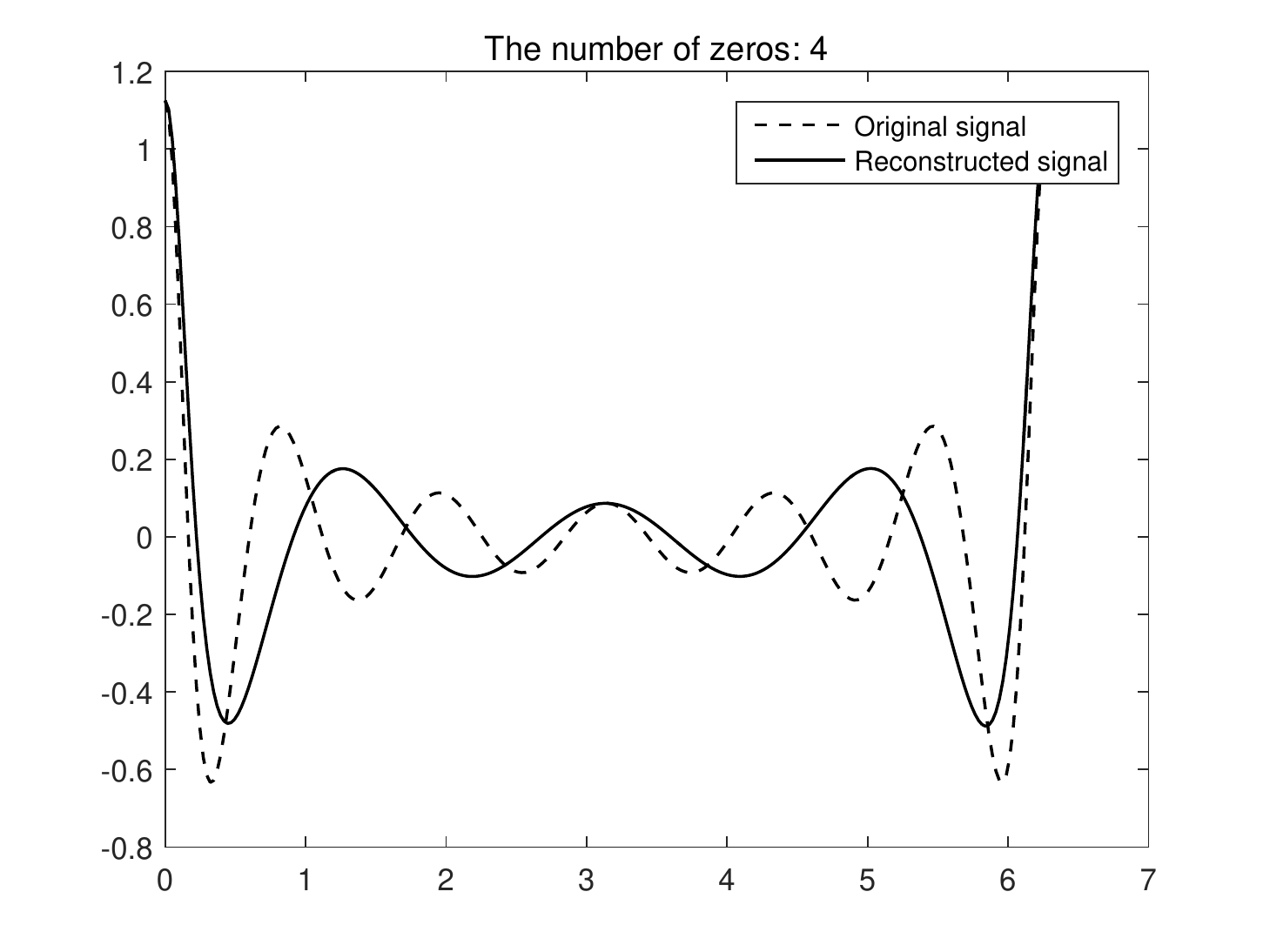}
\end{minipage}
\hfill
\begin{minipage}{0.3\linewidth}
  \includegraphics[width=2.2in,height=1.75in]{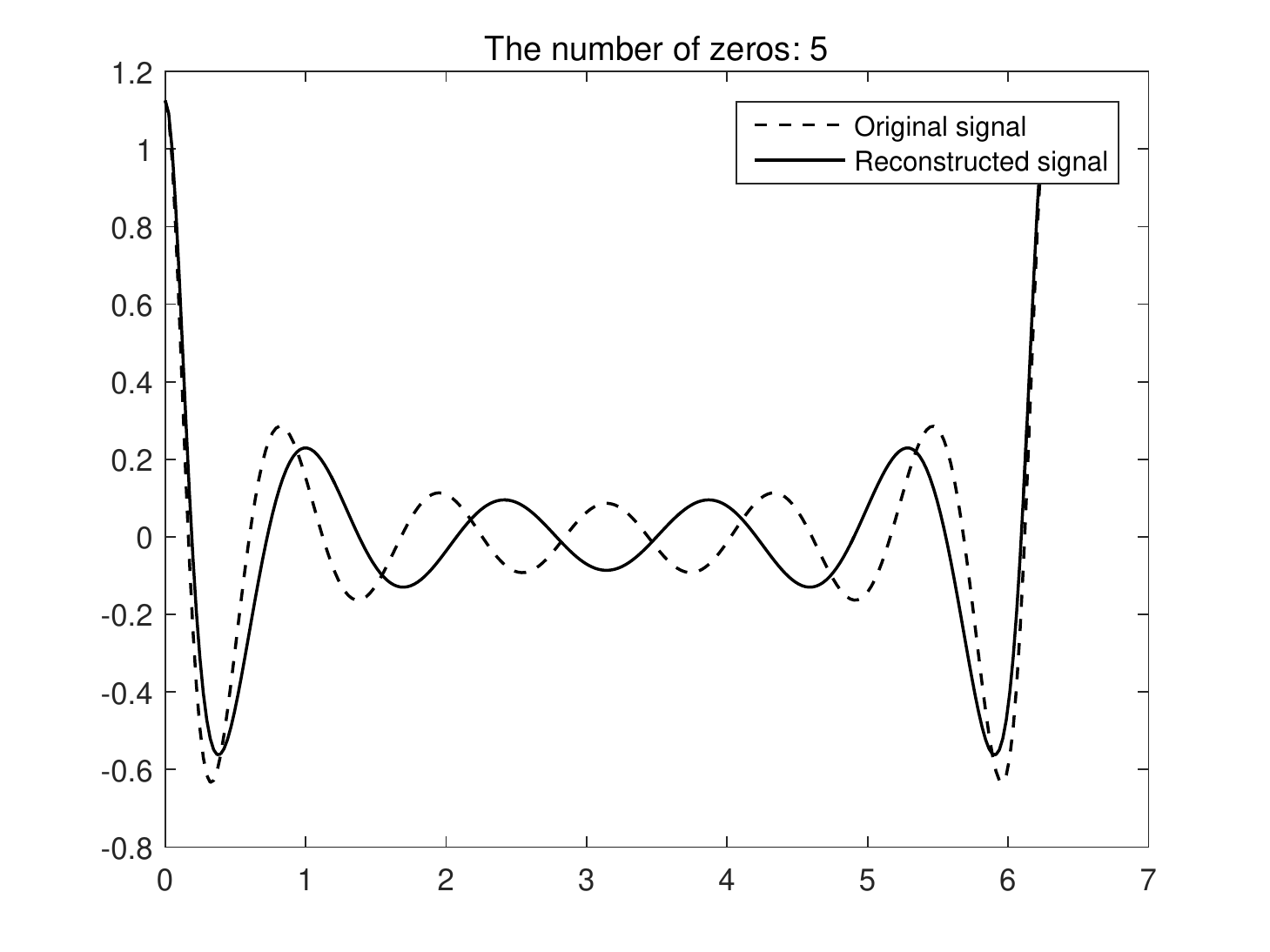}
\end{minipage}
\hfill
\begin{minipage}{0.3\linewidth}
  \includegraphics[width=2.2in,height=1.75in]{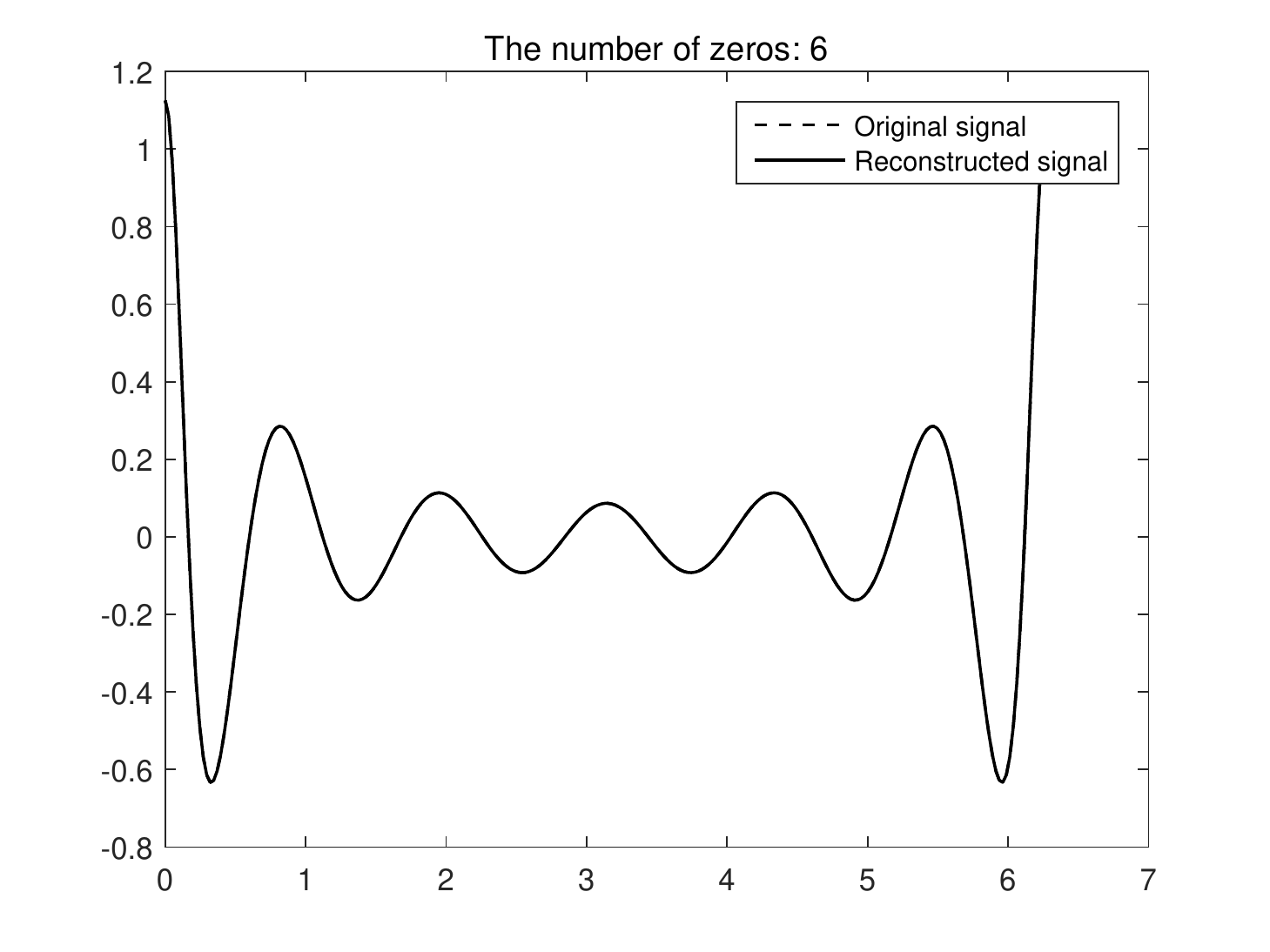}
\end{minipage}
\caption{256 sampling points}
\label{64}
\end{figure}

\begin{figure}[H]
\begin{minipage}{0.4\linewidth}
  \includegraphics[width=2.5in,height=2in]{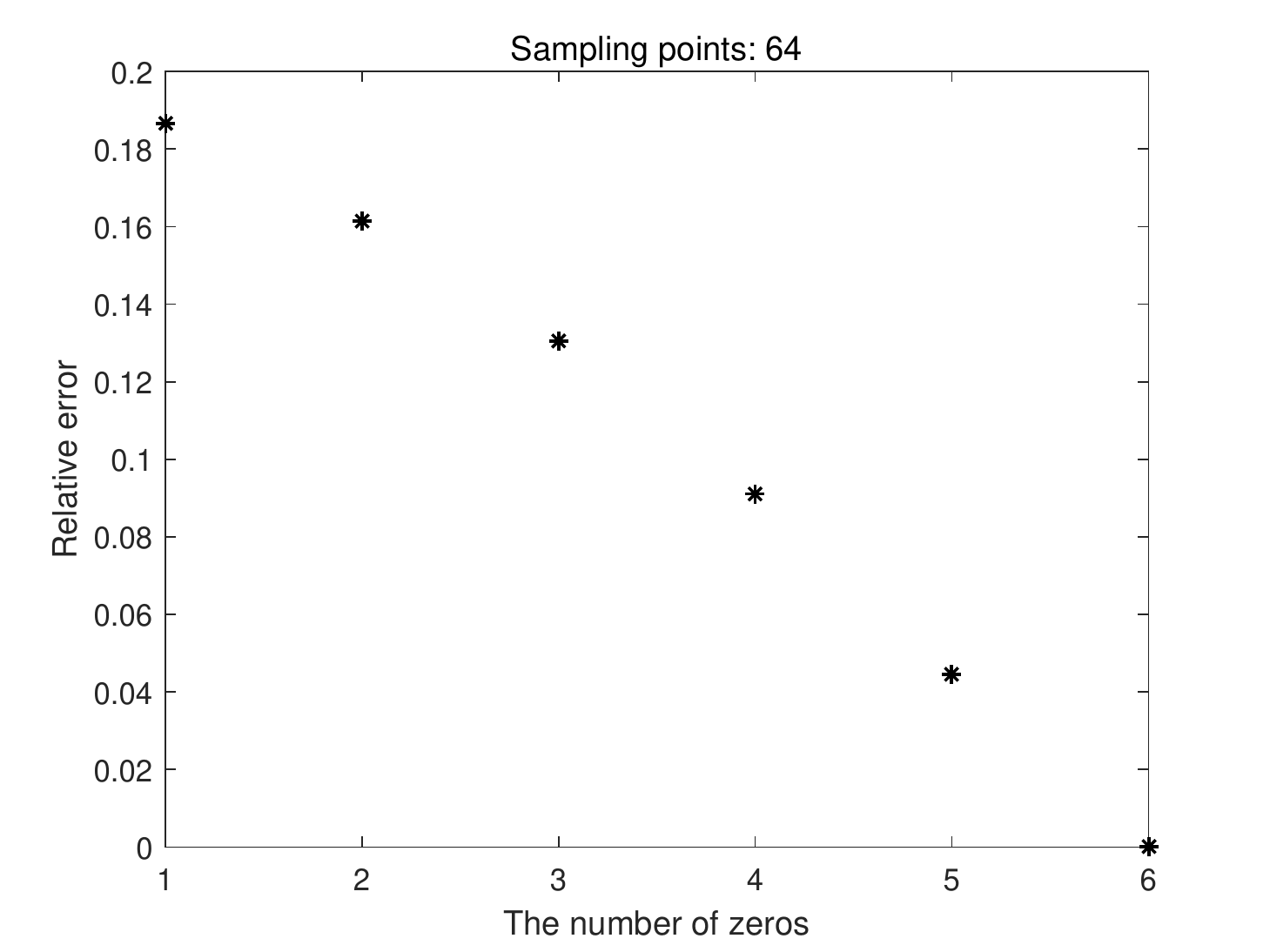}
\end{minipage}
\begin{minipage}{0.4\linewidth}
  \includegraphics[width=2.5in,height=2in]{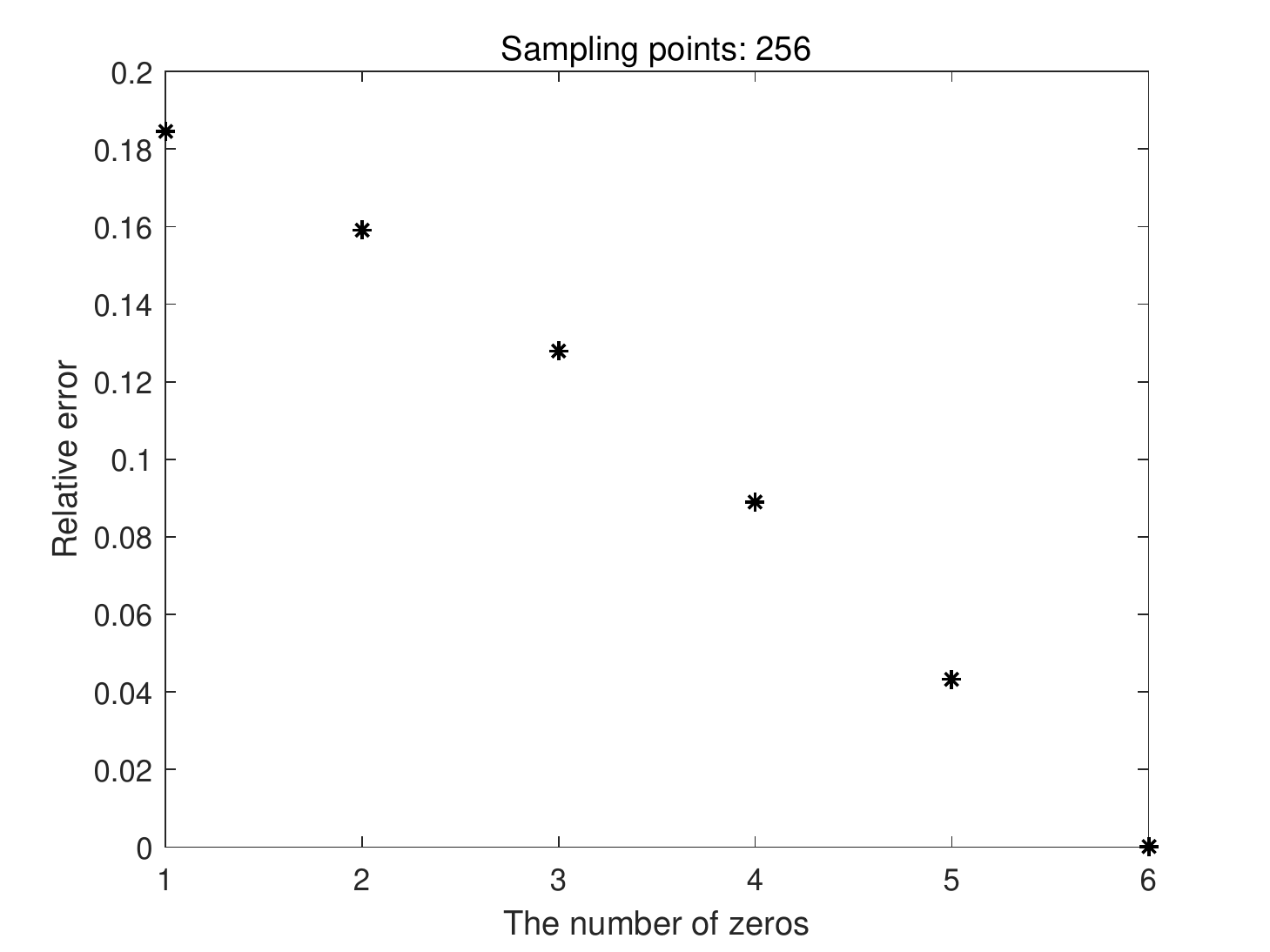}
\end{minipage}
\caption{Relative errors}
\label{0}
\end{figure}

\begin{table}[H]
    \centering
    \tiny
    \begin{tabular}{|c|c|c|c|c|c|c|}
    \hline
    \diagbox[width=14em,trim=l]{\tiny Sampling points}{\tiny The number of zeros} & 1 & 2 & 3 & 4 & 5 & 6 \\ \hline
    64 & 0.1866 & 0.1614 & 0.1305 & 0.0912 & 0.0446 & 1.8965$\times$ $10^{-5}$ \\ \hline
    256 & 0.1846 & 0.1591 & 0.1280 & 0.0899 & 0.0432 & 1.8667$\times$ $10^{-7}$ \\
    \hline
    \end{tabular}
    \caption{Relative errors}
    \label{error table}
\end{table}

\begin{example}
$$f(z)=\prod_{k=1}^{10}\frac{z-a_k}{1-\overline{a}_kz}$$
\end{example}
\begin{figure}[H]
\begin{minipage}{0.3\linewidth}
  \includegraphics[width=2.2in,height=1.75in]{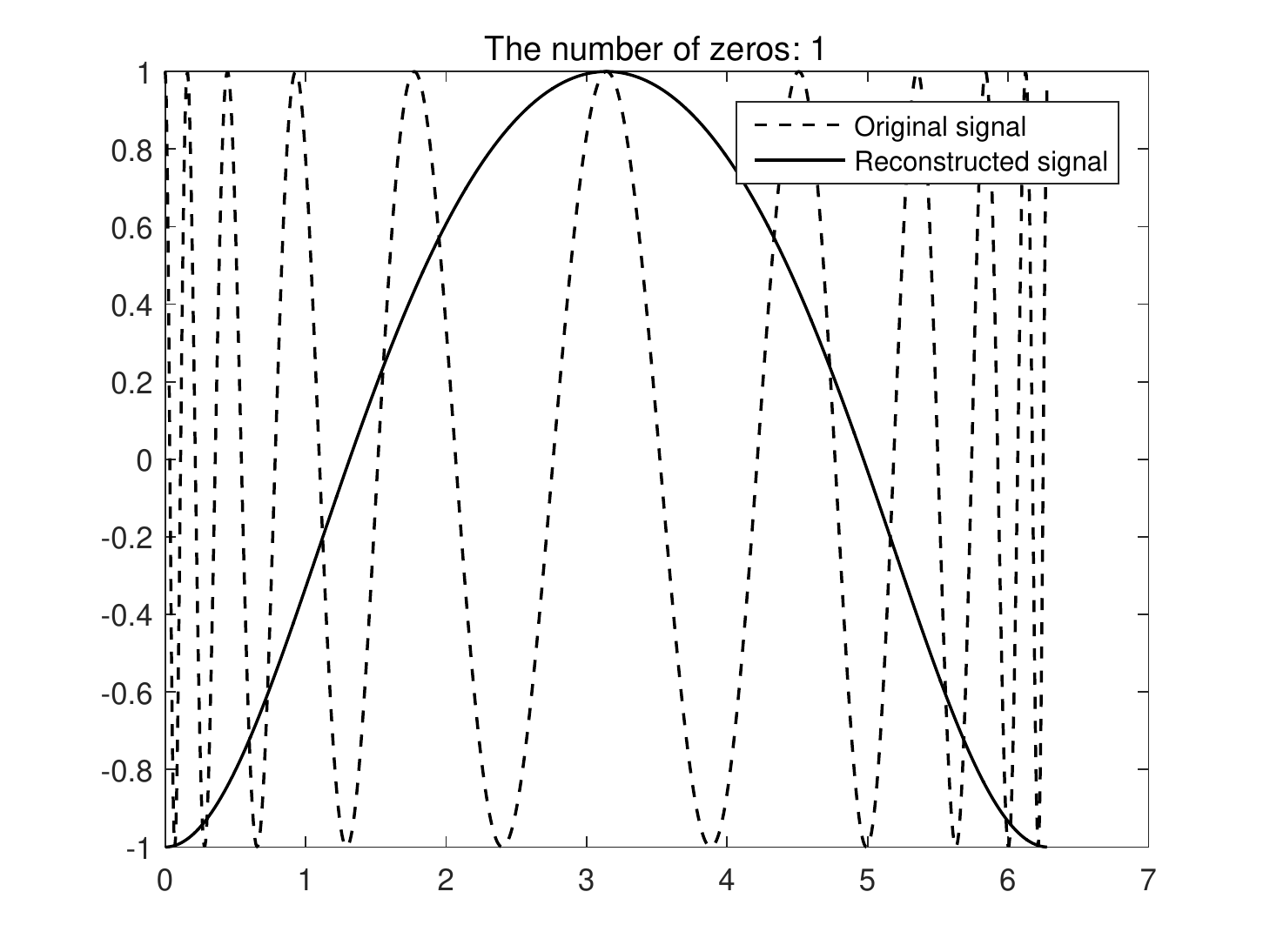}
\end{minipage}
\hfill
\begin{minipage}{0.3\linewidth}
  \includegraphics[width=2.2in,height=1.75in]{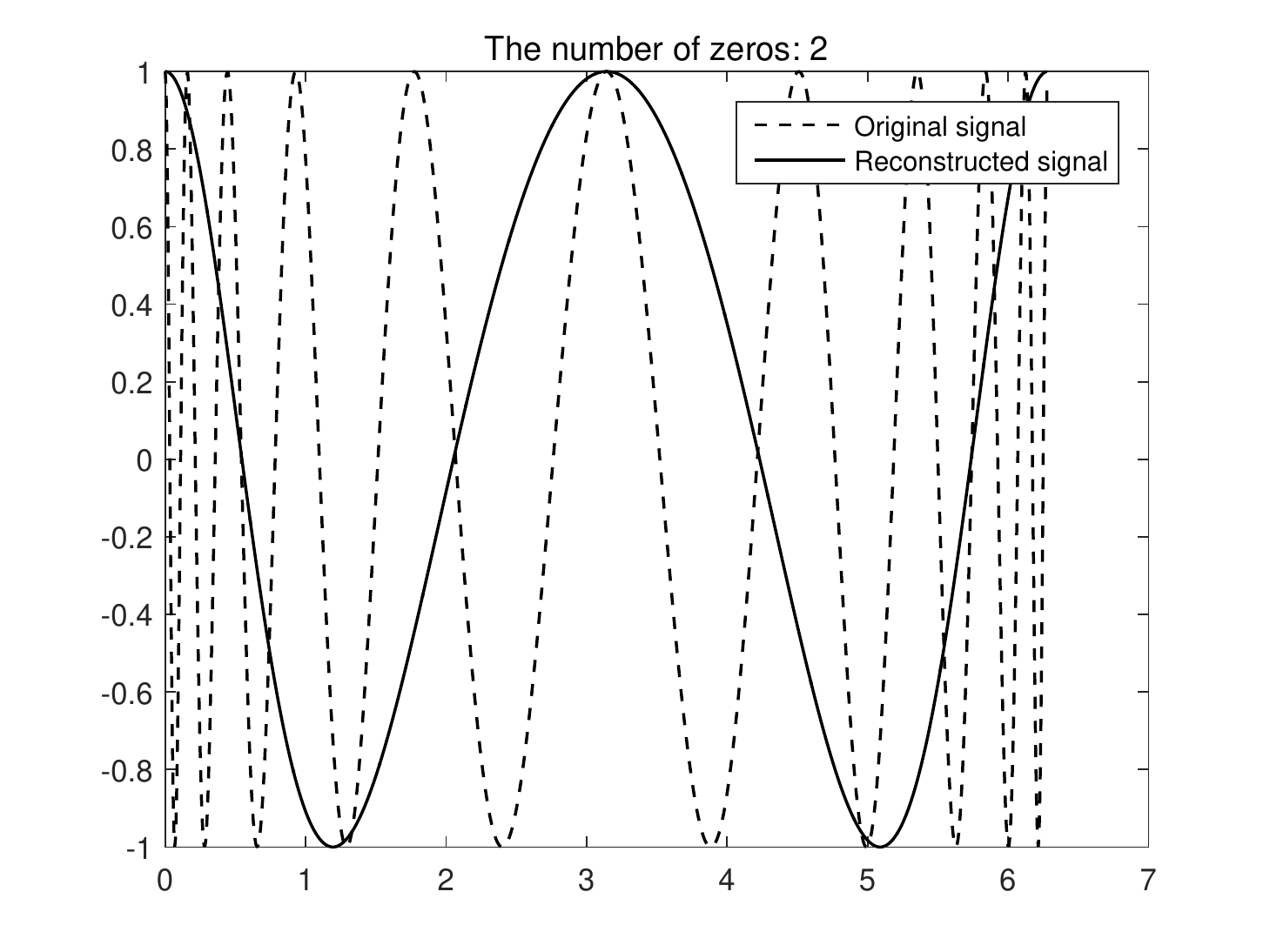}
\end{minipage}
\hfill
\begin{minipage}{0.3\linewidth}
  \includegraphics[width=2.2in,height=1.75in]{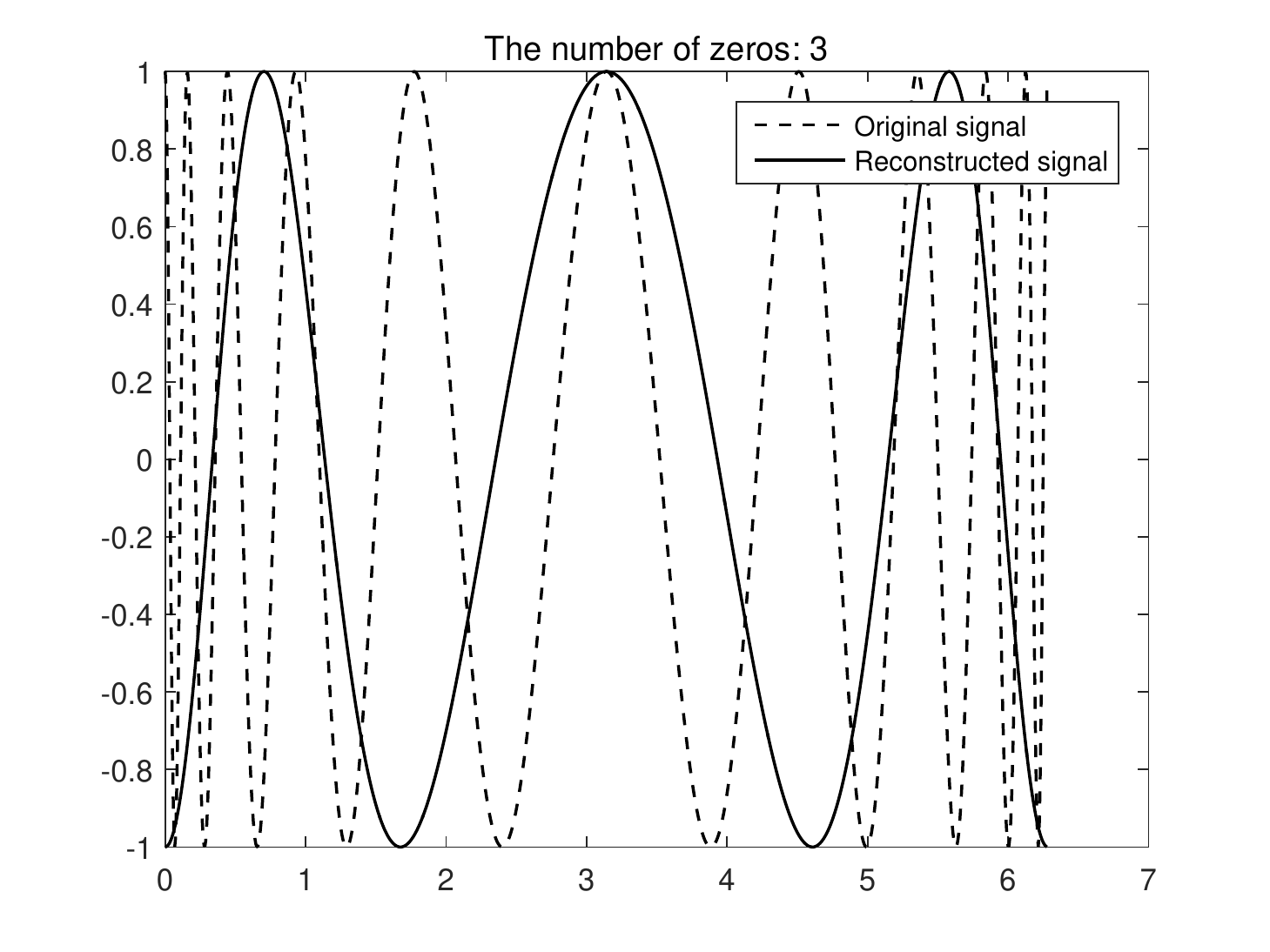}
\end{minipage}
\vfill
\begin{minipage}{0.3\linewidth}
  \includegraphics[width=2.2in,height=1.75in]{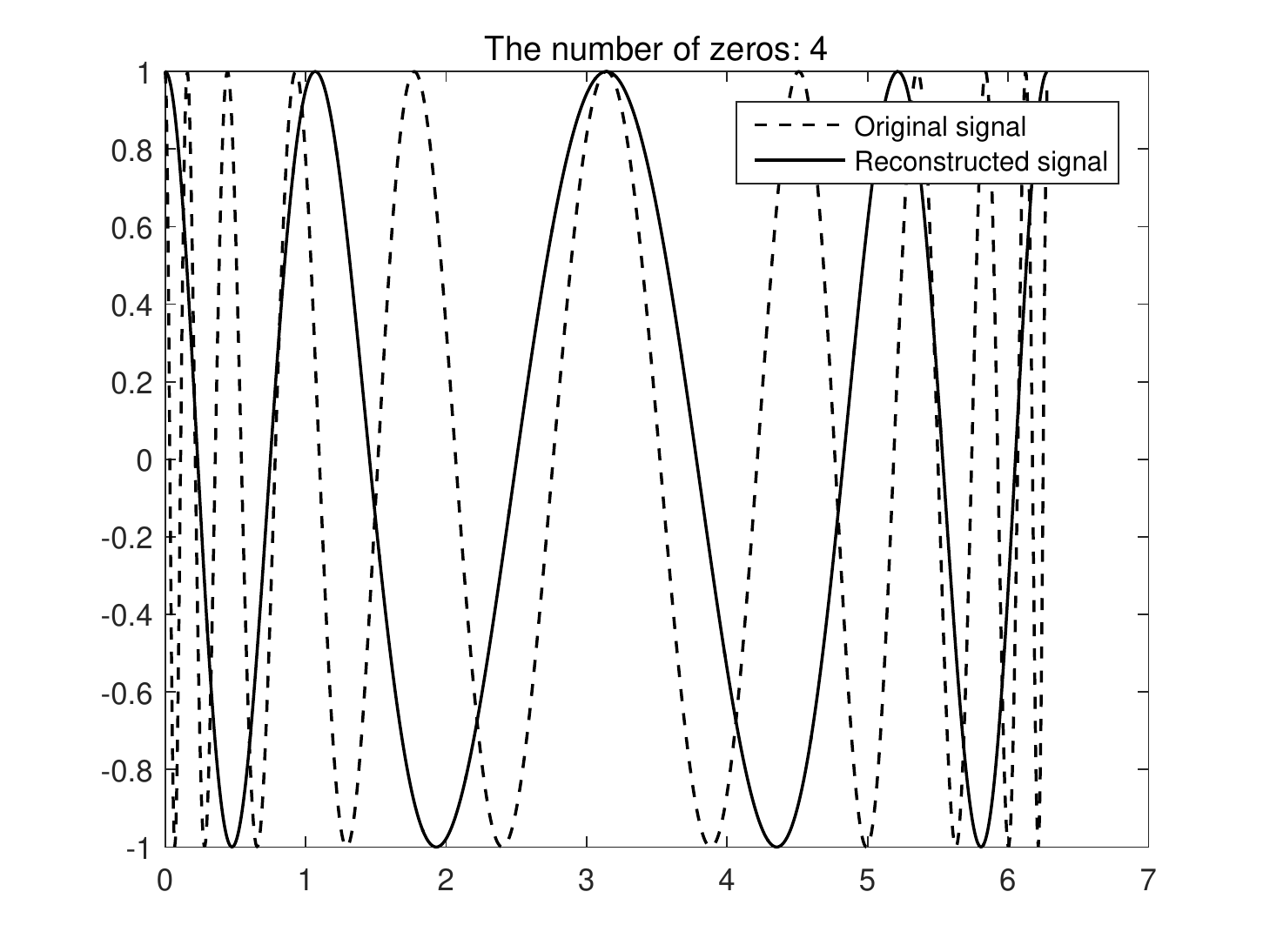}
\end{minipage}
\hfill
\begin{minipage}{0.3\linewidth}
  \includegraphics[width=2.2in,height=1.75in]{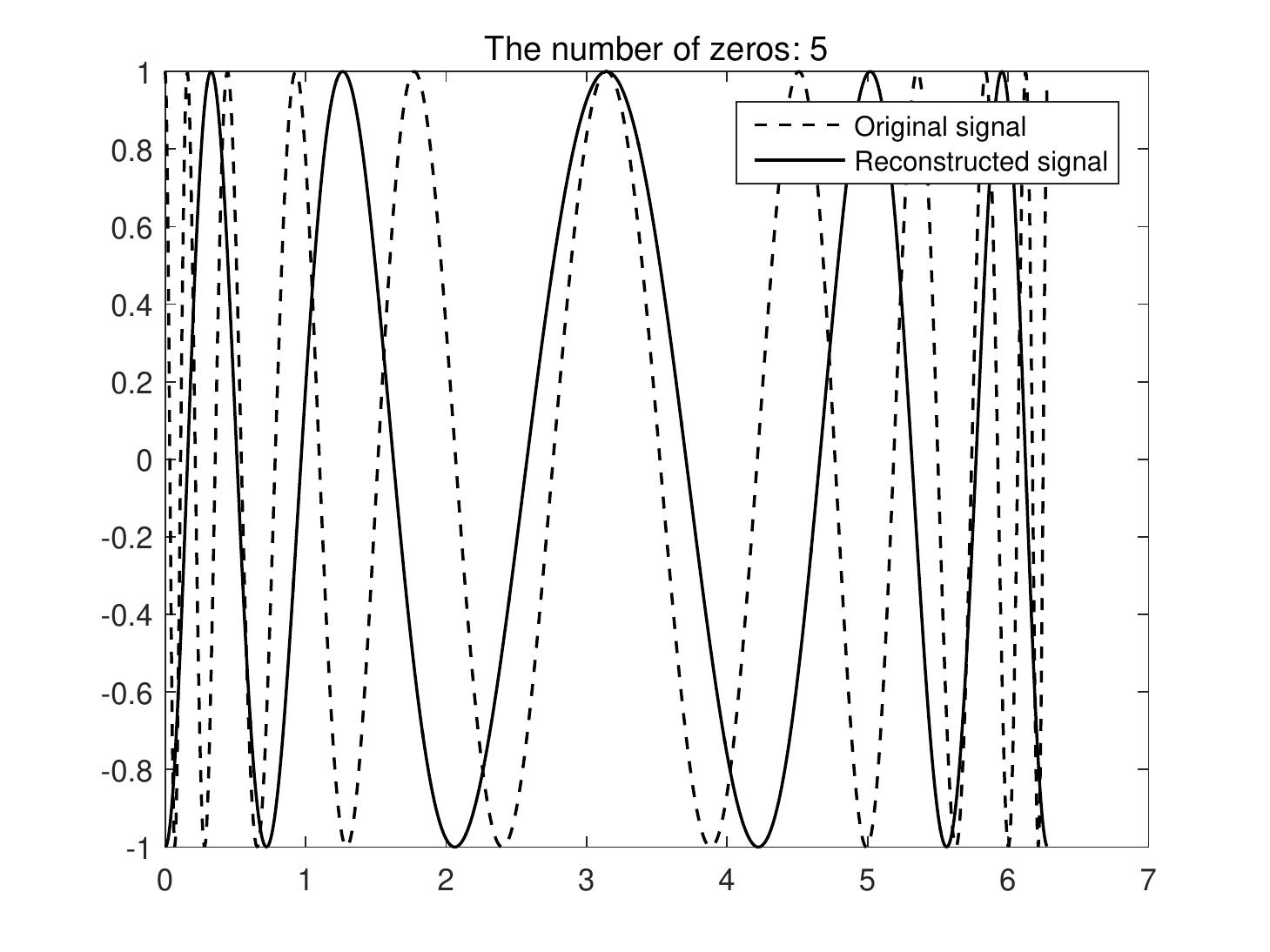}
\end{minipage}
\hfill
\begin{minipage}{0.3\linewidth}
  \includegraphics[width=2.2in,height=1.75in]{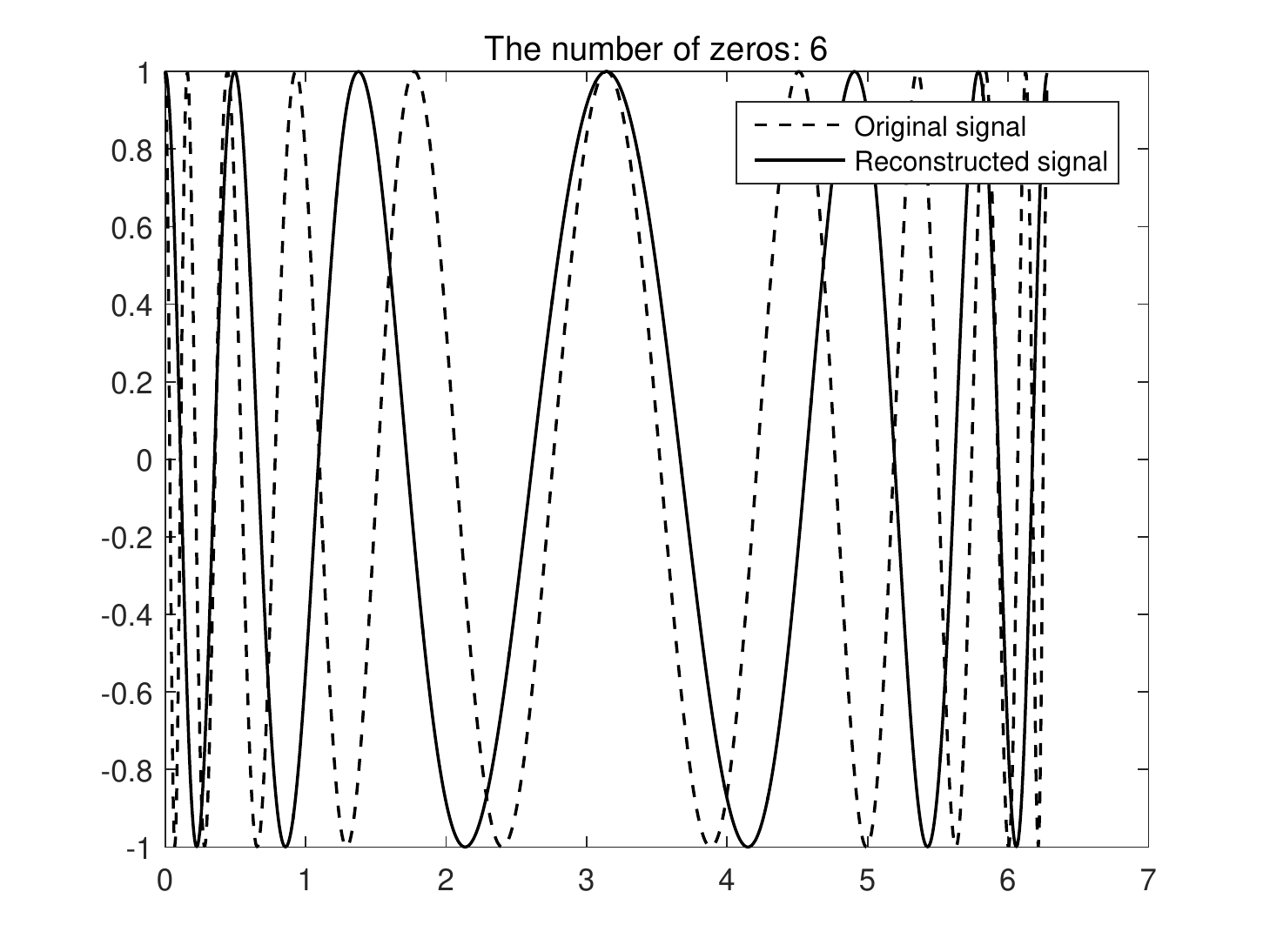}
\end{minipage}
\hfill
\begin{minipage}{0.3\linewidth}
  \includegraphics[width=2.2in,height=1.75in]{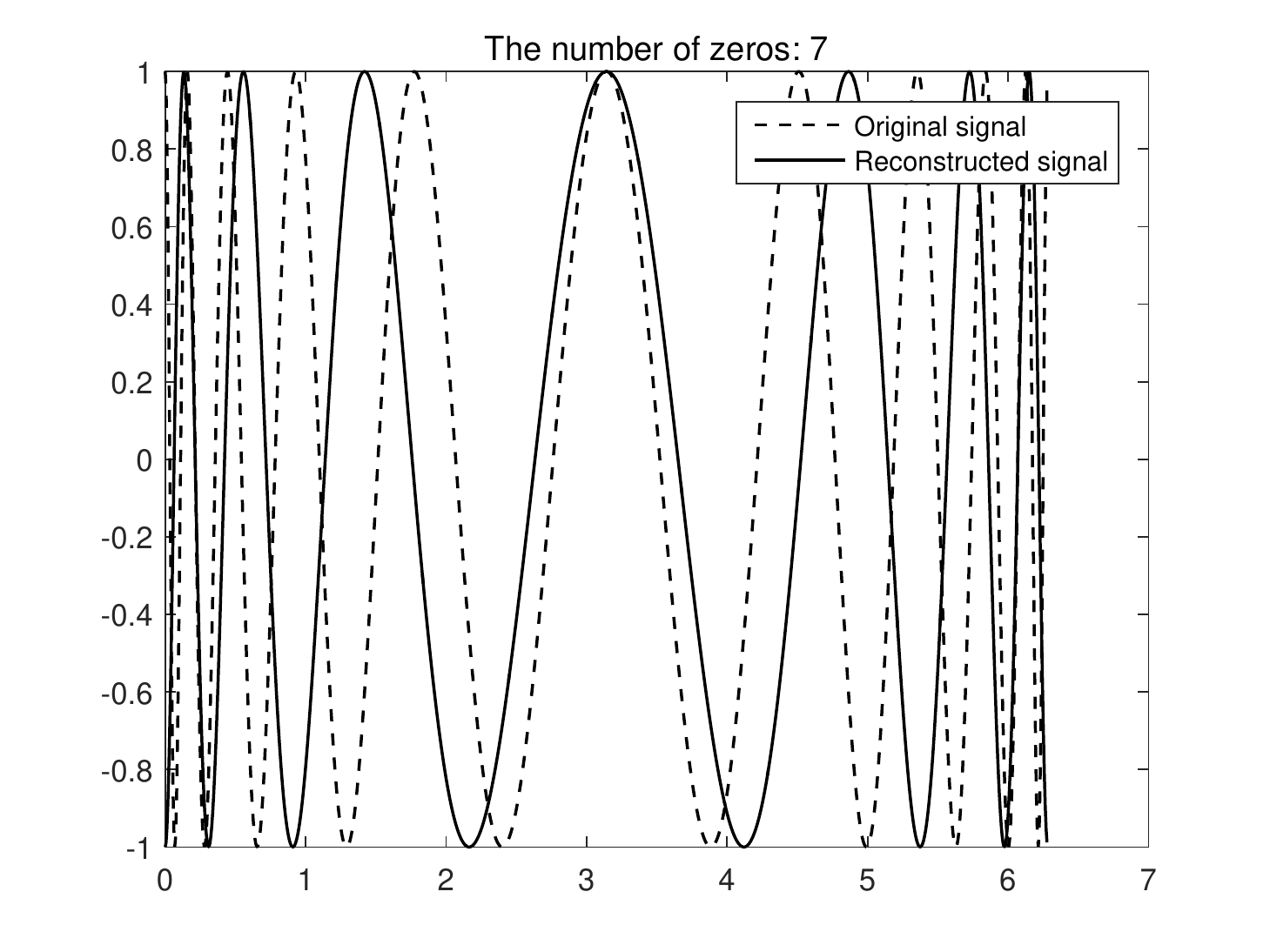}
\end{minipage}
\hfill
\begin{minipage}{0.3\linewidth}
  \includegraphics[width=2.2in,height=1.75in]{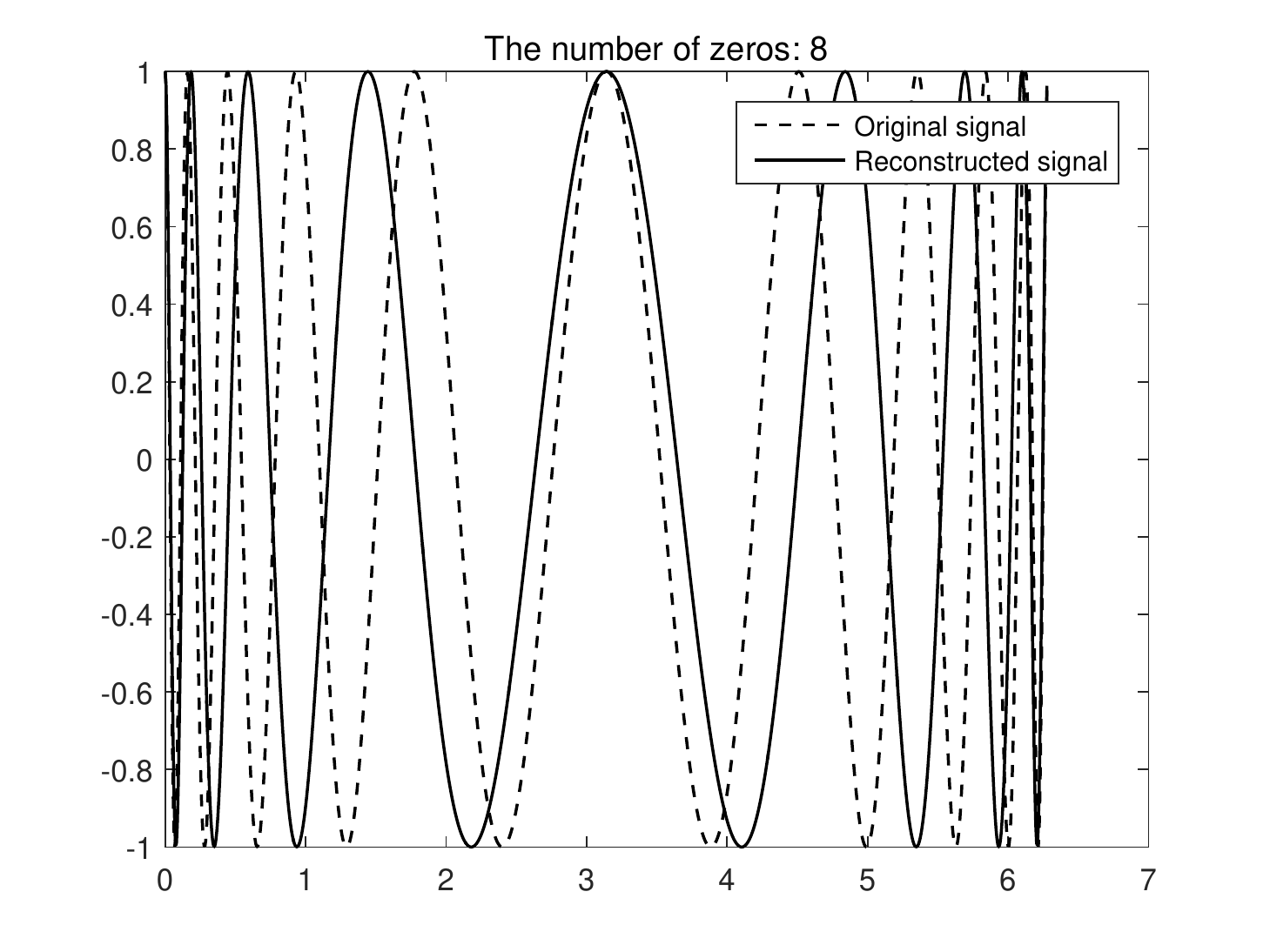}
\end{minipage}
\hfill
\begin{minipage}{0.3\linewidth}
  \includegraphics[width=2.2in,height=1.75in]{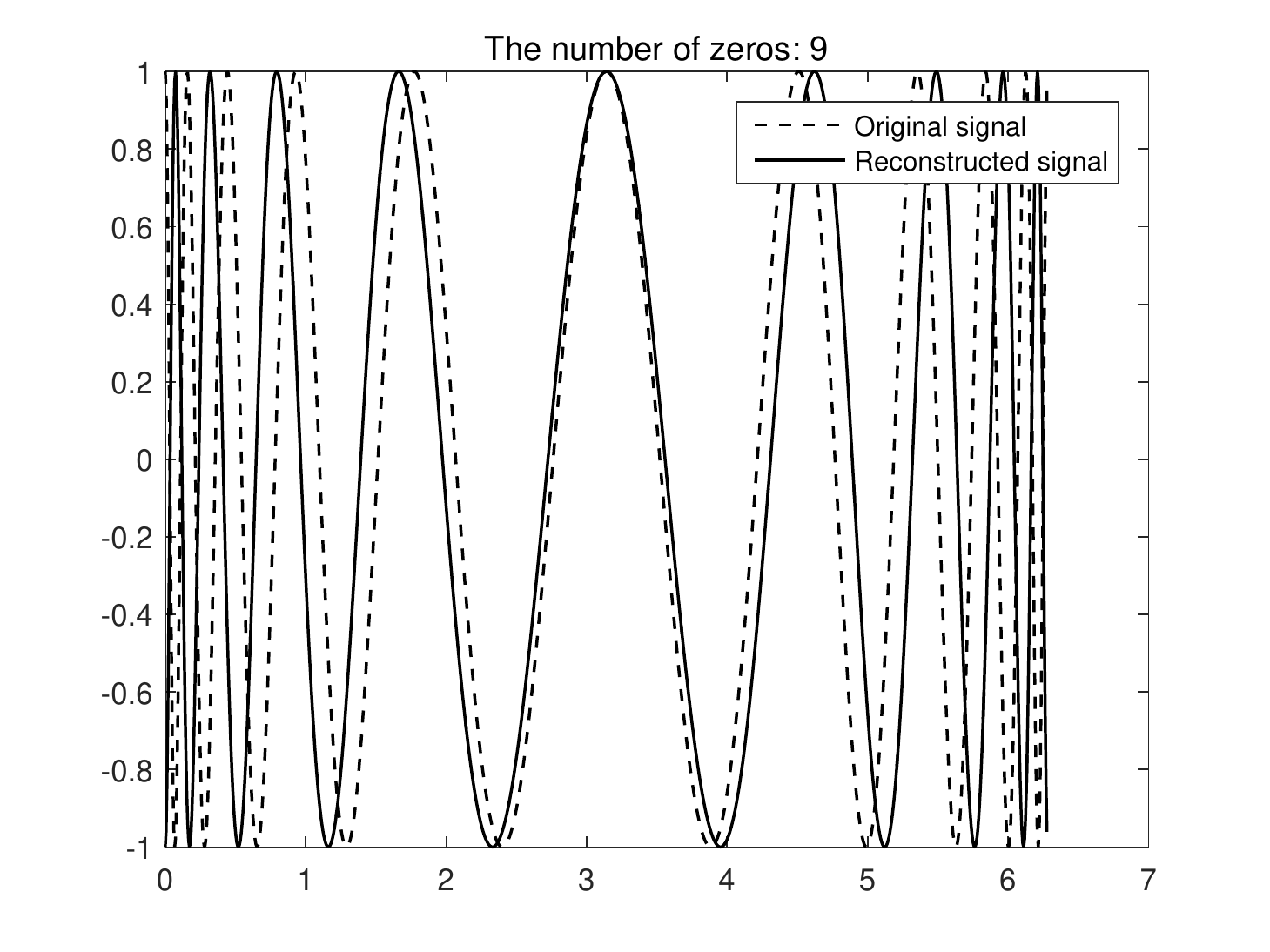}
\end{minipage}
\hfill

\begin{minipage}{0.45\linewidth}
  \includegraphics[width=2.2in,height=1.75in]{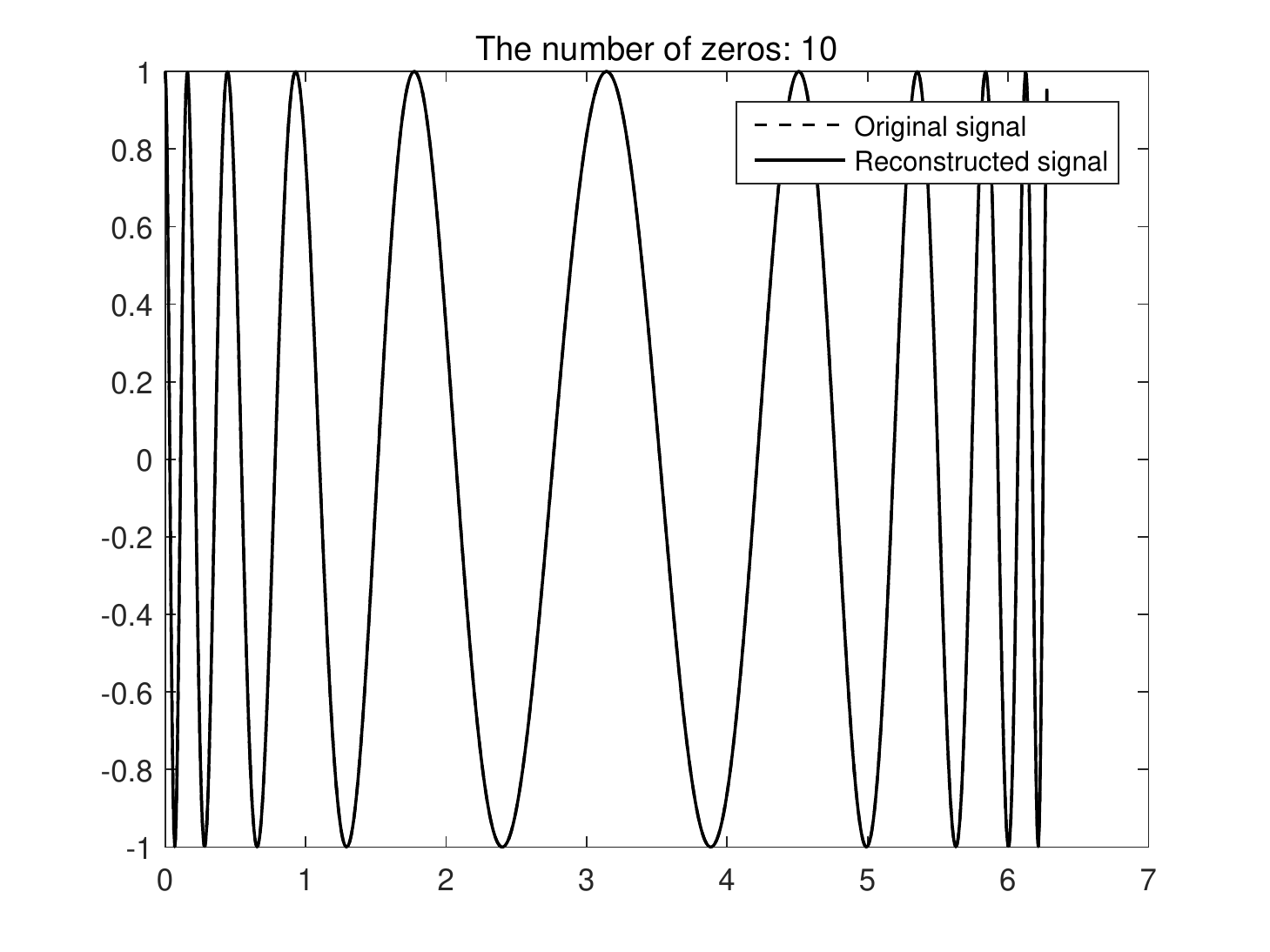}
\end{minipage}
\hfill
\begin{minipage}{0.45\linewidth}
  \includegraphics[width=2.2in,height=1.75in]{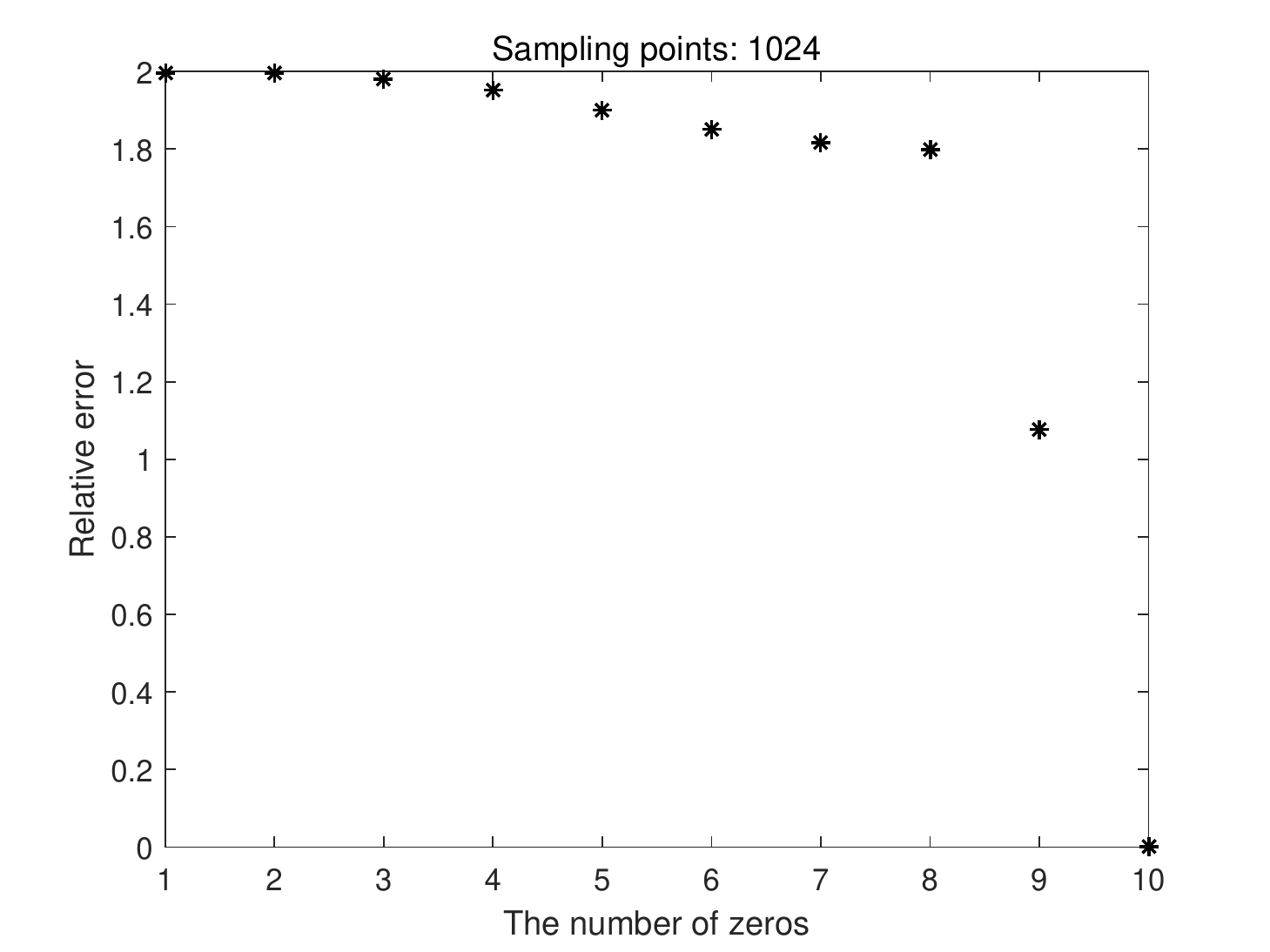}
\end{minipage}
\caption{1024 sampling points}
\label{1024}
\end{figure}

\end{document}